%% file: cohonil.tex
\documentclass[12pt]{amsart}

\usepackage{times,amsfonts,amsmath,amstext,amsbsy,amssymb,
amsopn,amsthm,upref}
\usepackage[T1]{fontenc}

\include{env_def}

\include{symb_def}
\def\dom{\text{dom}}
\begin{document}

\title%
{On the Cohomological equation for nilflows}

\author{Livio Flaminio} \author{Giovanni Forni}

\address{Math\'ematiques\\
  Universit\'e de Lille 1 (USTL)\\
  F59655 Villeneuve d'Asq CEDEX\\
  FRANCE}

\address{Department of Mathematics\\
  University of Toronto\\
  Toronto, ON  M5S 2E4 Canada}

\email
  {livio.flaminio@math.univ-lille1.fr}
\email
    {forni@math.toronto.edu}

\date{November 4, 2005}
    
\begin{abstract}
  \begin{sloppypar}
    Let $X$ be a vector field on a compact connected
    manifold $M$.  An important question in dynamical
    systems is to know when a function $g:M\to \R$ is a
    coboundary for the flow generated by $X$, i.e. when
    there exists a function $f: M\to \R$ such that $Xf=g$.
    In this article we investigate this question for
    nilflows on nilmanifolds. We show that there exists
    countably many independent Schwartz distributions $D_n$
    such that any sufficiently smooth function $g$ is a
    coboundary iff it belongs to the kernel of all the
    distributions $D_n$.
  \end{sloppypar}
 \end{abstract}

 \maketitle


\section{Introduction}%
\label{sec:Intoduction}

\subsection{The problem}
\label{ssec:intro1}

For a detailed discussion of the cohomological
equation for flows and tranformations in ergodic theory, we
refer the reader to \cite{MR2008435}. Here we limit
ourselves to a brief, self-contained introduction.

Let $X$ be a smooth vector field on a connected compact
manifold $\Cal M$ and let $(\phi_X^t)_{t\in \R}$ be the flow
generated by~$X$ on~$\Cal M$.  Many problems in dynamics and
ergodic theory (see loc. cit.) can be reconducted to the
study of the {\em cohomological equation}
\begin{equation} 
  \label{eq:cohoeq}
  Xu = f\,,
\end{equation}
i.e. the problem of finding a function $u$ on~$\Cal M$ whose
Lie derivative $Xu$ in the direction of~$X$ is a given
function~$f$ on~$\Cal M$.  Clearly a continuous solution $u$
of the equation \pref{eq:cohoeq} is only determined up to
function constant along the orbits. In addition, given the
value of $u$ at one point $p\in \Cal M$, the solution $u$ is
uniquely determined on the whole orbit of $p$ under the flow
$(\phi_X^t)_{t\in \R}$. In fact, the recurrence properties
of flow $(\phi_X^t)_{t\in \R}$ may as well forbid the
existence of measurable solutions to~\pref{eq:cohoeq}.
Hence the subtlety of the problem lies entirely in the fact
that, given $f$ in a certain regularity class, we may or may
not have solutions~$u$ in some other regularity class.

For example, observe that if the cohomological
equation~\pref{eq:cohoeq} admits a continuous solution $u$
then
\begin{equation}
  \label{eq:necess}
\mu(f)=0\text{ for all $X$-invariant Borel measures }\mu\in C^*(\Cal M).
\end{equation}
Thus invariant measures are obstructions to the existence of solutions
of~\pref{eq:cohoeq} in the continuous class. Livshitz's celebrated
theorem states that if $X$ is an Anosov vector field and $f$ is
H\"older continuos of exponent $\alpha\in (0,1)$, then the
necessary condition \pref{eq:necess} is also sufficient and that, in
this case, the solution $u$ is actually H\"older continuos of 
exponent $\alpha$. (This theorem has been generalized to other
hyperbolic dyanamical systems and to smoother classes of functions, see
for example \cite{MR1326374}).

For uniquely ergodic flows, i.e. for flows admitting a unique
invariant measure, one may wonder whether the above condition
\pref{eq:necess} is sufficient. It is well known that, for $X$ a
constant vector field on a torus, 
\begin{itemize}
\item if $X$ is Diophantine, for any $C^\infty$ function $f$ of
  average zero there exists a $C^\infty$ solution $u$ of
  equation~\pref{eq:cohoeq};
\item if $X$ is Liouvillean, then there exists some $C^\infty$
  function $f$ of average zero there for which there is no measurable
  weak\footnote{that is $(u,X\phi) = -(f,\phi) $ for all $\phi\in
    C^{\infty}(\Cal M)$} solution $u$ of equation~\pref{eq:cohoeq}.
\end{itemize}

In general, all invariant distributions (in the sense of
Schwartz), which are not necessarily signed measures, are
obstructions to the existence of smooth solutions of the
cohomological equation.

A distribution $D\in \mathcal D'(\Cal M)$ is called
$X$-invariant if $XD=0$ in the distributional sense, that
is, if $D(X\phi)=0$ for all $\phi \in C^{\infty}(\Cal M)$.
Let $\mathcal I^{\infty}_X(\Cal M)$ be the space of all
$X$-invariant distributions in $\mathcal D'(\Cal M)$.

By definition, if $u\in C^{\infty}(\Cal M)$ is a solution
of~\pref{eq:cohoeq}, we must have $D(f)=0$ for all
$D\in\mathcal I^\infty_X(\Cal M)$. By the same token, if $u$
is a solution of~\pref{eq:cohoeq} which is
$(\alpha+1)$-times differentiable, then we must have
$D(f)=0$ for all $X$-invariant distributions $D\in \mathcal
D'(\Cal M)$ of order $\alpha>0$.

Let $W^\alpha(\Cal M)$ be the Sobolev space of
square-integrable functions with square-integrable (weak)
derivatives up to order $\alpha>0$ and let $\mathcal
I^\alpha_X(\Cal M)\subset \mathcal I^{\infty}_X(\Cal M)$ be
the space of $X$-invariant distributions of Sobolev order
$\alpha>0$, i.e. which extend continuously to the Sobolev
space $W^{\alpha}(\Cal M)$.

If $u\in W^{\alpha+1}(\Cal M)$ is a solution
of~\pref{eq:cohoeq}, then $D(f)=0$ for all $D\in \mathcal
I^\alpha_X(\Cal M)$. It is then natural to ask if and when
invariant distributions are a complete set of obstructions
to the existence of differentiable solutions, i.e. when the
condition
\begin{equation}
  \label{eq:necess2}
  f\in W^\alpha(\Cal M) \quad\text{and}\quad D(f)=0\quad 
  \text{ for all } D\in \mathcal I^\alpha_X(\Cal M) \,,
\end{equation} 
is sufficient to guarantee the existence of a solution $u$
satisfying some regularity properties (for example $u \in
W^{\alpha-k}(\Cal M)$ for some $k\in \R^+$).

The main results of this paper is that this is the case for
all {\it Diophantine }nilflows on compact nilmanifolds.

\subsection{Main result}
\label{ssec:result}

Let~$ \nil $ be a $k$-step nilpotent real Lie algebra ($k\ge
2$) with a minimal set of generators $ \Cal E:= \{E_{1},
\dots , E_{n}\}\subset \nil$.  Let $\nil_{j}$, $j=1, \dots,
k$, denote the {\it descending central series} of $\nil$:
\begin{equation}
  \label{eq:descser}
  \nil_{1}=\nil,  \nil_{2}=[ \nil,\nil], \dots, \nil_{j}=[
  \nil_{j-1},\nil], \dots,  \nil_{k}\subset Z(\nil)\,,
\end{equation}
where $Z(\nil)$ is the center of $\nil$. 

Let $ \Nil$ be the connected and simply connected nilpotent
Lie group with Lie algebra $ \nil $. The corresponding Lie
subgroups $\Nil_j=\exp \nil_j=[\Nil_{j-1},\Nil]$ form the
descending central series of $\Nil$.

Let $\Gamma$ be a lattice in $\Nil$. It exists if and only
if $\Nil$ admits rational structure constants (see, for
example, \cite{MR0507234, CG:nilrep}).

A (compact) {\em nilmanifold} is a by definition a quotient manifold
$\manif:=\Gamma\backslash \Nil$ with $\Nil$ a nilpotent Lie group and 
$\Gamma\subset \Nil$ a lattice.

On a nilmanifold $\manif=\Gamma\backslash \Nil$, the group $
\Nil $ acts on the right transitively by right
multiplication. By definition, the {\it nilflow
  $(\phi_X^t)_{t \in \R}$ generated by $X\in \nil$ }is the
flow obtained by restriction of this action to the
one-parameter subgroups $(\exp t X)_{t \in \R}$ of $ \Nil $:
\begin{equation}
  \label{eq:nilflow}
  \phi_W^t(\Gamma x ) = \Gamma x \exp( t X ). 
\end{equation}
It is plain that nilflows on $\Gamma\backslash \Nil$
preserve the probability measure on $\Gamma\backslash \Nil$
given locally by the Haar measure. To simplify the notation,
the vector field on $\Gamma\backslash \Nil$ generating the
flow $(\phi_X^t)_{t \in \R}$ will also be indicated by $X$.

Every nilmanifold is a fiber bundle over a torus. In fact, the
group~$\overline{\Nil}=\Nil/[\Nil,\Nil]$ is abelian, connected and
simply connected, hence isomorphic to~$\R^n$ and~$\overline{ \Gamma}=
\Gamma/[\Gamma,\Gamma]$ is a lattice in~$\overline{\Nil}$. Thus we
have a natural projection
\begin{equation}
  \label{eq:projection}
  p: \Gamma\backslash \Nil\to 
  \overline{\Gamma}\backslash \overline{\Nil}  
\end{equation}
over a torus of dimension~$n$.

We recall the following:

\begin{theorem}[\cite{MR23:A3800}, \cite{MR29:4841}]
  \label{th:Green} The following properties are equivalent.
  \begin{enumerate}
  \item The nilflow $ ( (\phi_X^t)_{t \in \R}, \mu) $ on
    $\Gamma\backslash \Nil$ is ergodic.
  \item The nilflow $ (\phi_X^t)_{t \in \R} $ on $\Gamma\backslash
    \Nil$ is uniquely ergodic.
  \item The nilflow $ (\phi_X^t)_{t \in \R} $ on $\Gamma\backslash
    \Nil$ is minimal.
  \item The projected flow $ (\psi_{\bar X}^t)_{t \in \R} $ on $
    \overline{\Gamma}\backslash \overline{\Nil}\approx \T^n$ is an
    irrational linear flow on~$ \T^n $, hence it is (uniquely) ergodic
    and minimal.
  \end{enumerate}
\end{theorem}

The ``irrationality'' condition above in Theorem
\ref{th:Green}, $(4)$, refers to the rational structure
determined by the lattice $\overline{\Gamma}$.  Namely, if
the generators $\{E_1,\dots,E_n\}$ of~$\nil$ are chosen so
that the elements $\{\exp E_1,\dots,\exp E_n\}$ project onto
generators $\{\exp \overline{E}_1,\dots,\exp
\overline{E}_n\}$ of $\overline{\Gamma}$, then there exists
a vector $\Omega_X:=\left(\omega_1(X),
  \dots,\omega_n(X)\right)\in \R^n$ such that
\begin{equation}
\label{eq:barX}
\bar X =\omega_1(X) \overline{E}_1+\dots+\omega_n(X) \overline{E}_n
\end{equation}
and the condition means that $\omega_1(X), \dots,\omega_n(X)$ are linearly
independent over $\Q$. We shall call such an element $X\in \nil$ {\em
  irrational (with respect to $\Gamma$)}. 
  
  We have:
\begin{theorem} 
  \label{th:invardistrib} 
  Let $\manif=\Gamma\backslash \Nil$ be any compact
  nilmanifold obtained as a quotient of a $ k $-step
  nilpotent ($k\ge 2$) connected, simply connected Lie group
  $N$ by a lattice $\Gamma\subset N$. For any irrational
  $X\in \nil$, the space of $X$-invariant distributions
  $\mathcal I_X^\infty (\manif)$ has countable dimension.
  In fact, $\mathcal I_X^\infty (\manif)$ admits a countable
  basis $\Cal B_X$ of invariant distributions of Sobolev
  order $1/2$, in the sense that $\Cal B_X \subset \mathcal
  I^\alpha_X(\manif)$ for any $\alpha >1/2$.
\end{theorem}

We shall say that $X\in\nil $ is {\em Diophantine (with
  respect to $ \Gamma $) of exponent~$ \tau \geq 0 $ } if
the projection $\bar X$ of $X$ in $\nil/[\nil,\nil]$ is
Diophantine of exponent $ \tau\geq 0 $ for the lattice $\bar
\Gamma$ in the standard sense; namely, if there exists a
constant $ K>0 $ such that, for all $ M:=(m_{1},\dots,
m_{n}) \in \Z^{n}\setminus \{0\} $,
\begin{equation}
  \label{eq:DCintro} \vert \<M,\Omega_{X}\> \vert := \vert \sum_{i=1}^
  {n} m_{i} \omega_{i}(X) \vert \,\, \geq \,\, \frac{K}{ \vert
    M\vert ^{n-1+\tau}}\,\,.
\end{equation}
The subset of Diophantine elements $ X\in \nil $ of exponent
$ \tau\geq 0 $ will be denoted by $ \text{DC}_{\tau}(\nil)$.
The subset $ \text{DC}(\nil):= \cup_{\tau}
\text{DC}_{\tau}(\nil)$ of all Diophantine elements has full
Lebesgue measure in the Lie algebra $\nil$.

\smallskip Let $W_0^\alpha(\manif,X):= W^\alpha(\manif) \cap
\ker \mathcal I^\alpha_X(\manif)$ for any $\alpha>1/2$.  We
then have:
\begin{theorem}
  \label{th:main} 
  Let $\manif$, $\Gamma$, $N$ be as above. If $X\in
  \text{DC}_{\tau}(\nil)$ the following holds.
  \begin{enumerate}
  \item If $ \alpha>n + \tau-1/2$ and $ \beta<-1/2 $, there
    exists a linear bounded operator $ G_{X} :
    W^{\alpha}(\manif) \to W^{\beta}(\manif) $ such that,
    for all $ f\in W^{\alpha}(\manif) $, the distribution $
    u=G_{X}f $ is a weak solution of cohomological
    equation~{\rm \pref{eq:cohoeq}}.
  \item If $ \alpha>n+\tau $ and $ \beta<[\alpha-(n+\tau)]
    [(n+\tau)k+1]^ {-1} $, there exists a linear bounded
    operator $ G_{X} : W_0^{\alpha}(\manif,X) \to
    W^{\beta}(\manif) $, such that, for all $ f\in
    W^{\alpha}(\manif) $, the function $ u=G_{X}f $ is a
    solution of cohomological equation~{\rm
      \pref{eq:cohoeq}}.
  \end{enumerate}
\end{theorem}

\subsection{Motivations and applications}
\label{sec:motiv}

We have already mentioned that the cohomogical problem for a
flow (or a transformation) enters in many problems in
dynamical systems (KAM, time-changes, etc.).  We recall the
following definitions:
\begin{definition}
  Let $\Cal M$ be a compact connected manifold. A smooth
  vector field $X$ on $\Cal M$ is
  \begin{itemize}
  \item[$(a)$] {\em globally hypoelliptic (GH)}, if $Xu\in
  C^\infty(\Cal M)$ implies $u\in C^\infty(\Cal M)$ for any
  distribution $u\in \mathcal D'(\Cal M)$;
\item[$(b)$] {\em cohomology free (CF)}, if for all $f\in
  C^\infty(\Cal M)$ there exists a constant $c(f)\in \C$ and
  $u\in C^\infty(\Cal M)$ such that
  $$ Xu = f -c(f). $$
\end{itemize}
\end{definition}

In \cite{MR0320502} Greenfield and Wallach showed that, if a
vector field $X$ on $\Cal M$ is globally hypoelliptic, then
there exists an invariant volume form $\omega$ and the space
of $X$-invariant distributions is reduced to the line $\C
\omega$ (in particular the $X$-flow is uniquely ergodic).
The same is true, by a simple exercise, if $X$ is cohomology
free. In fact, if $X$ is cohomology free, then it is
globally hypoelliptic.

\begin{question}
  Is a globally hypoelliptic vector field on a compact connected
  manifold also cohomology free?
\end{question}

\begin{definition}[\cite{MR1858535,MR2008435}]
  A vector field $X$ on a compact connected manifold $\Cal
  M$ is said {\em $C^\infty$-stable} if the subspace $\{ Xu
  \mid u \in C^\infty(\Cal M)\}$ is closed in $C^\infty(\Cal
  M)$.
\end{definition}

It is plain that for $C^\infty$-stable vector fields the
answer to the above question is positive. Absence of
stability is related to fast approximation by periodic
systems: the major example of $ C^\infty$-{\em un}stable are
Liouvillean constant vector fields on tori (loc.cit.). In
\cite{MR0320502}, Greenfield and Wallach showed that if a
Killing vector field of a Riemannian metric on $\Cal M$ is
globally hypoelliptic, then the manifold $\Cal M$ is a torus
and that on the flat $2$-torus any Killing globally
hypoelliptic vector field must be a constant Diophantine
vector field. Motivated by these results they conjectured:

\begin{conjecture}[\cite{MR0320502}]
  If a compact, connected manifold $\Cal M$ admits a
  globally hypoelliptic vector field $X$ then $\Cal M$ is a
  torus and $X$ is smoothly conjugate to a constant
  Diophantine vector field.
\end{conjecture}

In fact, if $\Cal M$ is a torus, by \cite{MR1737551} any
globally hypoelliptic vector field $X$ is smoothly conjugate
to a constant Diophantine vector field.

A related conjecture due to A. Katok is the following:
\begin{conjecture}[\cite{MR1858535,MR805843,MR2008435}] 
  If a compact, connected manifold $\Cal M$ admits a
  cohomology free vector field $X$ then $\Cal M$ is a torus
  and $X$ is smoothly conjugate to a constant Diophantine
  vector field.
\end{conjecture}

In support of Katok's conjecture Federico and Jana Rodriguez
Hertz have recently proved in \cite{RHRH:coho} that any
(compact, connected) manifold $\Cal M$ admitting a
cohomology free vector field fibres over a torus of
dimension equal to the first Betti number of $\Cal M$.

Theorems \ref{th:invardistrib} and \ref{th:main} immediately
imply that, within the class of nilflows, both conjectures
are true. As remarked by the Rodriguez Hertz's in
\cite{RHRH:coho}, nilflows are the simplest class of flows
for which their theorem yields no non-trivial information.

\begin{sloppypar}
  Another motivation for the study of cohomological
  equations and invariant distributions comes from
  ``renormalization dynamics'' (see \cite{f:icm}). In fact,
  whenever we can define, for a family of ``parabolic
  flows'', an effective renormalization dynamics on a
  suitable moduli space, the analysis of the action of the
  renormalization on the bundle of invariant distributions
  over the moduli space may allow to determine the
  asymptotics of the ergodic averages for the typical flow
  in the given family.  This program was carried out in
  \cite{MR1888794} for conservative flows on (higher genus)
  surfaces (related to interval exchange transformations),
  in \cite{MR2003124} for horocycle flows on hyperbolic
  surfaces and in \cite{flafor:heis} for nilflows on the
  3-dimensional Heisenberg group. These results have shown
  that the phenomenon of polynomial deviations of ergodic
  averages for interval exchange transformations, discovered
  by A. Zorich \cite{zorich1}, \cite{zorich2},
  \cite{zorich3} and rigorously proved in \cite{MR1888794},
  is shared by other fundamental examples of parabolic
  flows. In principle, it should be possible to carry out
  this program by purely dynamical methods based on
  renormalization. In fact, in the case of interval exchange
  transformations, Marmi, Moussa and Yoccoz \cite{MMY} have
  been able to replace the harmonic analysis methods applied
  by the second author in his study of the cohomological
  equation \cite{MR99d:58102} by a renormalization approach
  (based on the Rauzy-Veech-Zorich induction). However, no
  renormalization dynamics is currently available for
  general nilflows.
\end{sloppypar} 

\section{Irreducible Unitary Representations}%
\label{sec:irrepsofnil}

\subsection{Kirillov's Theory}
\label{ssec:kirillov}

Let~$ \Nil $ be the connected and simply connected nil\-potent Lie
group with Lie algebra~$ \nil $.  By Kirillov theory, all the
irreducible unitary representation of~$ \Nil $ are parametrized
by the {\em coadjoint orbits}~$ \Cal O\subset \nil^{\ast} $, i.e.
by the orbits of the {\em coadjoint action} of~$ \Nil $
on~$ \nil^\ast $ defined by
\begin{equation}
  \label{eq:coadact}
  \Ad^*(g)\lambda =
  \lambda\circ \Ad(g^{-1})\qquad g\in \Nil, \lambda\in \nil^\ast.
\end{equation}
For~$ \lambda\in \nil^* $, the skew-symmetric bilinear form
\begin{equation}
  \label{eq:skewform}
  B_\lambda(X,Y)=\lambda([X,Y])
\end{equation}
has a {\em radical~$ \mathfrak r_\lambda $} which coincides
with the Lie subalgebra of the subgroup of ~$ \Nil $
stabilizing~$ \lambda $; thus the form~\pref{eq:skewform} is
non-degenerate on the tangent space to the orbit~$ \Cal
O\subset \nil^{\ast} $ of $ \lambda $ and defines a
symplectic form on~$ \Cal O $.

A {\em polarizing (or maximal subordinate) subalgebra for $
  \lambda $} is a maximal isotropic subspace $ \m \subset \nil $
for the form $ B_\lambda $ which is also a subalgebra~of~$ \nil
$. In particular any polarizing subalgebra for the linear form~$
\lambda $ contains the radical~$ {\mathfrak r}_{\lambda}$.  If
$\m$ is a polarizing subalgebra for the linear form~$ \lambda $,
the map
$$
\exp T \mapsto \exp 2 \pi \imath \lambda(T),\qquad T\in \m, 
$$
yields a one-dimensional representation, which we denote by $\exp
2 \pi \imath \lambda$, of the subgroup $M=\exp \m\subset \Nil$.

To a pair $ \Lambda:=(\lambda,\m) $ formed by a linear
form~$ \lambda \in \nil^\ast $ and a polarizing subalgebra~$\m$
for~$ \lambda $, we associate the unitary representation
\begin{equation}
  \label{eq:induced}
\pi_{\Lambda} =\text{Ind}_{\exp \m}^N (\exp 2 \pi \imath \lambda).  
\end{equation}
These unitary representations are irreducible; up to
unitary equivalence, all unitary irreducible representations of~$
\Nil $ are obtained in this way.  Furthermore, two pairs $
\Lambda:=(\lambda,\m) $ and $ \Lambda':=(\lambda',\m') $ yield
unitarily equivalent representations $ \pi_{\Lambda} $ and $
\pi_{\Lambda'} $ if and only if $\lambda$~and~$\lambda'$ belong
to the same coadjoint orbit~$ \Cal O\subset \nil^{\ast} $.

The unitary equivalence class of the representations of
the group~$ \Nil $ determined by the coadjoint orbit~$ \Cal O $
will be denoted by~$\Pi_{\Cal O} $, while we set
$$
\Pi_{\lambda} =\{\pi_\Lambda \mid \Lambda=(\lambda,\m), \text
{ with $\m$ polarizing subalgebra for $\lambda$} \} .
$$

\begin{definition} 
  Let $\Nil$ be a connected, simply connected Lie group of
  Lie algebra $\nil$.  The \/ {\em derived representation}
  $\pi_*$ of a unitary representation $\pi$ of $\Nil$ on a
  Hilbert space $H_\pi$ is the Lie algebra representation of
  $\nil$ on $H_\pi$ defined as follows.  For every $X\in
  \nil$,
\begin{equation}
  \pi_* (X)\, := \, 
  \text{\rm strong-}\lim_{t\to  0} (\pi(\exp tX ) - I)/t \,.
\end{equation}
\end{definition} 

It can be proved that the derived representation $\pi_*$ of
the Lie algebra $\nil$ on $H_\pi$ is {\it essentially
  skew-adoint } in following sense. For all $X\in \nil$, the
linear operators $ \pi_* (X)$ are essentially skew-adjoint
with common invariant core the subspace of $C^\infty(H_\pi)
\subset H_\pi$ of $C^{\infty}$-vectors in $H_\pi$. We recall
that a vector $v\in H_\pi$ is $C^{\infty}$ for the
representation $\pi$ if the function $g\in \Nil \mapsto \pi
(g) v\in H_{\pi}$ is of class $C^\infty$ as a function on
$\Nil$ (with values in a Hilbert space).

Suppose that $N$ is the semi-direct product $A\ltimes N'$ of a
normal subgroup $N'$ and an abelian group $A$ and that $\pi'$
is a unitary irreducible representation of $N'$ on an
Hilbert space $H'$; then the derived representation $\pi_*$
of the induced representation $\pi=\text{Ind}_{N'}^N (\pi')$
has a simple description in terms of $\pi'_*$: in fact, up
to unitary equivalence, $H_\pi\approx L^2(A,H')$; furthermore the
subgroup $A$ acts by translations, and hence for $f\in
L^2(A,H')$ we have
\begin{equation}
  \label{eq:semi1}
  \left(\pi_*(X) f\right)(a) = \frac{d} {dt}\left. f(a\exp
    tX)\right|_{t=0}, 
  \quad a\in A,\; X \in \mathfrak a := \text{Lie}(A);
\end{equation}
the infinitesimal action of $N'$ is a pointwise action given
explicitely by:
\begin{equation}
  \label{eq:semi2}
  \left(\pi_*(Y) f\right)(a) = \pi'_*\left(\Ad(a)Y\right) (f(a)), 
  \quad    a\in A,\;Y \in {\mathfrak n}' := \text{Lie}(N').
\end{equation}
Since, for any $a\in A$, the operator $\Ad(a)$ acting on $
{\mathfrak n}'$ is unipotent, the right-hand side in the formula
above is a polynomial function in the variable $a\in A$. In
fact, if $\ell$ is the degree of nilpotency of $\ad(X)$ and
$a=\exp X$, we have, for all $Y \in {\mathfrak n}' = \text{Lie}(N')$,
\begin{equation}
  \label{eq:semi3}
  \left(\pi_*(Y) f\right)(a) = 
  \sum_{j=0}^\ell \frac 1 {j!} \pi' _*\left( \ad(X)^{j}Y\right) (f(a))\,.
\end{equation}

\subsection{Coadjoint orbits of maximal rank} 

Let~$ \nil $ be a $ k $-step nilpotent real Lie algebra on
$n$~generators $E_{1}, \dots , E_{n}$.  Let $\nil_{j}$,
$j=1, \dots, k$, denote the {\it descending central series}
of $\nil$:
$$  
\nil_{1}=\nil, \nil_{2}=[ \nil,\nil], \dots,
\nil_{j}=[\nil_{j-1},\nil], \dots, \nil_{k}\subset Z(\nil).
$$
In this section we characterize coadjoint orbits that
correspond to unitary representations which do not factor
through the quotient $\Nil/\! \exp\nil_k$. Such coadjoint
orbits and the induced unitary representations will be
called of {\it maximal rank}. Since the Lie group $\Nil /\!
\exp\nil_k$ is $(k-1)$-step nilpotent, the analysis of unitary
representations of maximal rank is sufficient to treat by induction
all unitary representations of a given $k$-step nilpotent Lie group.

We shall make the fundamental assumption 
that the coadjoint orbit $\Cal O \subset \nil^{\ast}$ has
maximal rank (see below). Before stating this assumption, a few
lemmas and definitions.

Let
\begin{equation}
\label{eq:nilperp}
\begin{aligned}
\nil_{k-1}^{\perp}(\lambda)\,&=\,\{ X\in \nil \,\mid\,
B_{\lambda} (X,\nil_{k-1})=0\, \}, \\
\nil_{k-1}^{\perp}(\Cal O)\,&=\,\nil_{k-1}^{\perp}(\lambda)\,,
\quad \text{ for any } \lambda\in \Cal O \,.
\end{aligned}
\end{equation}

Since the restriction of $\lambda\in \nil^{\ast}$ to the
centre $Z(\nil) $ does not depend on the choice of the
linear form $ \lambda\in \Cal O $ and since $ [\nil,
\nil_{k-1}] = \nil_{k}\subset Z(\nil)$, the restriction $
B_{\lambda} \vert \nil\times \nil_{k-1} $ and the subspace $
\nil_{k-1} ^ {\perp}(\lambda) $ depend only on the coadjoint
orbit~$\Cal O $.  Consequently, the second definition
in~\pref{eq:nilperp} is well-posed.

\begin{lemma}
  \label{lem:trivial}
  Let $ \Cal O $ be a coadjoint orbit and $ \lambda\in \Cal O $.
  Then
  \begin{enumerate}
  \item $ \mathfrak r_\lambda\subset \nil_{k-1} ^{\perp}(\Cal O)
    $\,;
  \item $\nil_2 \subset \nil_{k-1} ^{\perp}(\Cal O)$\,;
  \item $\nil_{k-1} ^{\perp}(\Cal O)$ is a sub-algebra.
  \end{enumerate}
\end{lemma}
\begin{proof}
  \begin{sloppypar}
    Condition $(1)$ is immediate from the definitions.
    Conditions $(2)$ and $(3)$ follow from Jacobi's identity
    and the inclusion $\nil_{k}\subset Z(\nil)$.  In fact,
    if $B_{\lambda}(T_1, \nil_{k-1})= B_{\lambda}(T_2,
    \nil_{k-1}) =\{0\} $ then $B_{\lambda}([T_1,T_2],
    \nil_{k-1})= \lambda([[T_1,T_2], \nil_{k-1}])\subset -
    \lambda([[T_2,\nil_{k-1}], T_1]) -
    \lambda([[\nil_{k-1},T_1], T_2])= \{0\}$, since
    $[T_j,\nil_{k-1}]\subset \nil_k\subset Z(\nil)$.
  \end{sloppypar}
\end{proof}

\begin{lemma}
  \label{lem:nonmaxrk}
  Let $ \Cal O $ be a coadjoint orbit and $ \lambda\in \Cal O $.
  The following properties are equivalent:
  \begin{enumerate}
  \item the restriction $ \lambda \vert\nil_{k} $ is
    identically zero
  \item $\nil_{k-1} ^{\perp}(\Cal O) =\nil$
  \item the projection of $\nil_{k-1} ^{\perp}(\Cal O)$ on  $
    \nil/\nil_{2} $ is surjective
  \end{enumerate}
\end{lemma}
\begin{proof}
  (\,1$\Rightarrow$2\,)~\,Since $ [\nil, \nil_{k-1}] =
  \nil_{k}$, if the restriction $ \lambda \vert\nil_{k} $ is
  identically zero, by definition, we have $\nil_{k-1}
  ^{\perp}(\lambda) =\nil$. 

  (\,2$\Rightarrow$1\,)~\,From $ [\nil, \nil_{k-1}] =
  \nil_{k}$, if $\nil_{k-1} ^{\perp}(\lambda) =\nil$, we
  have $ \lambda \vert\nil_{k} =0$.

  (\,3$\Leftrightarrow$2\,)~\,Any subalgebra of $\nil$ containing
  $\nil_2$ and which projects onto $\nil/\nil_{2} $ must concide
  with the full algebra~$\nil$.
\end{proof}
If the conditions of the above lemma are satisfied, any~$ \pi\in
\Pi_{\Cal O} $ factors through a representation of~$ \Nil/\exp
\nil_k $, a $(k-1)$-step nilpotent group. Having in mind an
induction process on the degree of nilpotency $k$, we shall
assume that we are not in this case.
\begin{definition}
  If none of the conditions of the Lemma~\ref{lem:nonmaxrk} are
  satisfied we will say that the coadjoint orbit~$ \Cal O $, any
  linear form $ \lambda\in \Cal O $ and any irreducible
  representation $ \pi \in \Pi_{\Cal O} $ have {\em maximal rank
  }(equal to $ k $).
\end{definition}

Thus from now on we shall focus on coadjoint orbits~$ \Cal O
$ of maximal rank.

\subsection{Adapted representations}  
Let~$ \Cal O $ be any coadjoint orbit of maximal rank. In
Section~\ref{sec:cohom} we shall study the cohomological
equation for a ``generic'' admissible $X\in \nil$ restricted
to an irreducible unitary representation $\pi \in \Pi_{\Cal
  O}$.  There we shall see, in fact, that the set of
admissible $X\in \nil$ equals $\nil\setminus
\nil_{k-1}^\perp (\Cal O)$.
  
The following lemma produces, for any coadjoint orbit $\Cal
O$ of maximal rank and any $X\in \nil\setminus
\nil_{k-1}^\perp (\Cal O)$, a special irreducible unitary
representation $\pi\in \Pi_{\Cal O}$ adapted to our goal of
proving a priori estimates for the cohomological equation
restricted to that representation.

\begin{lemma}
  \label{lemma:cohomeqone} 
  Let $ X\in\nil\setminus \nil_{k-1}^\perp (\Cal O)$ and let
  $Y\in \nil_{k-1}$ be any element such that $
  B_{\lambda}(X,Y) \not = 0$ for all $\lambda\in \Cal O $.
  There exists a codimension 1 ideal $\nil'\subset \nil$
  with $ X\not\in \nil'$ and a unitary (irreducible)
  representation~$ \pi\in \Pi_{\Cal O} $ with the following
  properties:
  \begin{enumerate}
    \item the representation $\pi$ is obtained inducing from
    $\Nil':=\exp {\nil}'$ to $\Nil$ a unitary irreducible
    representation~$ \pi' $ of $\Nil'$ on a Hilbert space
    $H'$;
    \item the derived the representation $\pi_*$ of the Lie
    algebra $\nil$ satisfies \smallskip
  \begin{enumerate}
    \item
      $ \pi_*(X)= \frac{\partial}{\partial t} $ on~$ L^{2} (\R,
      H', dt) $,
    \item
      $ \pi_*(Y) =2\pi \imath \, B_{\lambda}(X,Y)\, t \,\I_{H'}
      $ on $ L^{2} (\R, H', dt) $;
  \end{enumerate}
   \smallskip
 \item there exists a constant $C>0$ such that 
  \begin{equation}
  \label{eq:Xproj}
  \vert \<X,U\> \vert \geq  C^{-1} \, 
  \frac{ \vert B_{\lambda}(X,Y)\vert }
  {  \Vert \lambda\vert \nil_k \Vert }\,.
 \end{equation}
 \end{enumerate}
 where $U\in \nil$ denote the unit normal vector to $\nil'$
 with respect to an euclidean product $\<\cdot, \cdot\>$ on
 $\nil$ fixed once for all.
\end{lemma}

\begin{proof}
  Since, by Lemma~\ref{lem:trivial} we have $\nil_2\subset
  \nil_{k-1}^\perp(\mathcal O)$, for any $\lambda\in\Cal O
  $, the subspace $ \nil' =\{ T\in \nil \mid B_{\lambda}(T,
  Y) =0\} $ is a codimension $ 1 $ ideal of $ \nil $
  depending only on the coadjoint orbit $ \Cal O $. Let
  $U\in \nil$ be a normal unit vector to $\nil'$ with
  respect to an euclidean product $\<\cdot, \cdot\>$ on
  $\nil$ fixed once for all. There is an orthogonal
  decomposition
  \begin{equation}
  X= \<X,U\> U \,+\,  W   \,,
  \end{equation}
  for some $W\in \nil'$. Since by definition of $\nil'$ we
  have $B_\lambda(W,Y)=0$, it follows immediately that
  \begin{equation}
    \label{eq:BXY}
    B_\lambda(X,Y)= \<X,U\> B_\lambda(U, Y) \,.
  \end{equation}

  There exists a constant $C>0$ such that $\Vert U\Vert
  =\Vert Y \Vert =1$ implies that $\Vert [U,Y] \Vert \leq C$
  and $Y\in\nil_{k-1}$ implies $[U,Y]\in \nil_k$.  Hence
  \begin{equation}
    \label{eq:BUY}
    \vert B_\lambda(U, Y) \vert =
    \vert \lambda ([U,Y])
    \vert \leq   C \,\Vert \lambda\vert \nil_k \Vert \,.
  \end{equation}
  The lower bound \pref{eq:Xproj} follows immediately from
  \pref{eq:BXY} and \pref{eq:BUY}.
 
  Let $ \m \subset \nil' $ be a polarizing subalgebra for
  the restriction $ \lambda'=\lambda\vert \nil' $ of~$
  \lambda$ to~$ \nil' $. Then $ \m $ is also a polarizing
  subalgebra for $ \lambda $ on $ \nil $ (\cite{CG:nilrep},
  Prop.~1.3.4), since, by definition, we have $ \mathfrak
  r_{\lambda} \subset \nil' $ for all $\lambda\in\Cal O $.
  
  Set
  $$ 
  \Lambda=(\lambda,\m),\quad \Lambda'=(\lambda',\m),\quad
  \pi =\pi_{\Lambda} \quad \text{ and~}\pi' =\pi_{\Lambda'}.
  $$
  By induction in stages (see for example \cite{CG:nilrep},
  \S 2.1), we have:
  $$
  \pi = \text{Ind}_{\exp \m}^\Nil \big(\exp ({\imath
    \lambda})\big)\approx \text{Ind}_{\exp \nil'}^\Nil
  \Big(\text{Ind}_{\exp \m}^{\exp \nil'} \big(\exp ({\imath
    \lambda'})\big) \Big)= \text{Ind}_{\exp \nil'}^\Nil
  (\pi').
  $$
  From \pref{eq:semi1}-\pref{eq:semi2}, if we identify
  with~$\R$ the one-parameter subgroup generated by~$X$ and
  if we denote by~$H'$ the Hilbert space on with
  $\pi_{\Lambda'}$ acts, up to a unitary equivalence, the
  representation $\pi_{\Lambda}$ coincides with the
  representation of~$\Nil$ on~$ L^{2} (\R, H', dt) $
  infinitesimally given by:
  $$
  \pi_*(X) f(t) = \tfrac{\partial}{\partial t} f(t),\qquad
  \pi_*(Y) f(t)= \pi'_*\big(\text{Ad}(\exp tX)Y\big)f(t), \qquad
  $$
  for $Y\in \nil'$ and all $f$ in a dense subspace of~$ L^{2}
  (\R, H', dt)$, such as the subspace $C^{\infty }_0 (\R,
  H')$ of compactly supported $C^\infty$ functions from~$\R$
  to~$H'$.
  
  Since~$[X,Y]\in \nil_{k}\subset Z(\nil)$ we have $ [X,[X,Y]]=0 $ and
  therefore 
  \begin{equation}
    \label{eq:adxz}
    \text{Ad}(\exp tX)Y = Y + t[X,Y];
  \end{equation}
  we also have~$\nil_{k}\subset \nil'$, hence
  $\nil_{k}\subset Z(\nil') \subset \m' $ and $\pi'|
  \nil_{k}=\exp \imath \lambda$; from \pref{eq:semi3}, we
  obtain that $ \pi_*(Y)$ is the operator pointwise given
  by:
  $$
    \pi_*(Y)= t\, \pi'_*([X,Y]) \,+\, \pi'_*(Y) = 2\pi\imath \,t
    \,B_{\lambda}(X,Y) \, \I_{H'} \, + \,\pi'_*(Y) \,.
  $$

  We claim that there exist $ \lambda\in \Cal O $ and a
  representation $\pi'\in \Pi_{\lambda'} $ such that~$ \pi'_*
  (Y)= 0 $. This will prove the lemma.

  Since $ [X,Y] \in Z(\nil) $, the value $ B_{\lambda} (X,Y) $
  does not depend on $ \lambda\in \Cal O $.  Let $ \lambda\in \Cal
  O$. From~\pref{eq:adxz} we have
  $$
  \big(\Ad^\ast(\exp t_0X)\lambda\big)(Y)=\lambda
  \big(\Ad(\exp(-t_0X))Y \big) =\lambda (Y )- t_0 B_\lambda(X,Y)
  $$
  which vanishes for $t_0= \lambda (Y )/B_\lambda(X,Y)$.
  Thus we can and shall assume that we have chosen a linear
  form~$\lambda\in \Cal O $ satisfying~$ \lambda(Y)= 0 $.
  
  Let $\nil_k''=\{ T\in \nil_k \mid \lambda(T)= 0 \}$.
  By hypothesis~$\nil_k''$ is a central ideal in~$\nil$ of
  codimension~1  in~$\nil_k$.  We set
  $$
  \bar\nil'= \nil'/\nil_k''
  $$
  and we denote by 
  $\bar\lambda'$ the linear
  form
  induced on 
  $\bar\nil_k'$ by 
  $\lambda'$. We
  also denote by 
  $$
  p':\Nil'=\exp\nil'\to \Nil'/\exp\nil_k''
  $$
  the canonical projection. We have $\Pi_{\lambda'}=\{
  \bar\pi'\circ p' \mid \bar\pi'\in \Pi_{\bar\lambda'}\}$ since
  the central ideal $\nil_k''$ is contained in every
  subalgebra~$\m\subset \nil'$ polarizing for~$\lambda'$ and every
  representation $\pi'\in\Pi_{\lambda'}$ is trivial on
  $\exp\nil_k''$. But now observe that the projection~$\bar Y$
  of~$Y$ in~$\bar\nil'$ belongs to the centre~$Z(\bar\nil')$
  of~$\bar\nil'$. In fact, if~$T\in \nil'$, we have $[T,Y]\in
  \nil_k$ and then the condition $B_\lambda(T,Y)=0$ is equivalent
  to $[T,Y] = 0 \mod \nil_k''$. We conclude that for every
  $\bar\pi'\in \Pi_{\bar\lambda'}$ we have $\bar\pi'_*(\bar Y) =
  \imath \bar \lambda'(\bar Y)=0$ and consequently $\pi'_*(Y) =
  0$ for any $\pi'\in \Pi_{\lambda'}$.

  This concludes the proof.
\end{proof}

\section{The cohomological equation in representation}
\label{sec:cohom}

\subsection{Distributions and Sobolev spaces}%
Let $\Nil$ be a connected, simply connected Lie group.  Let
$\pi $ be a unitary representation of $\Nil$ on a Hilbert
space $H_\pi$. Let $\pi_*$ denote the derived representation
of the Lie algebra $\nil$ of $\Nil$ on $H_\pi$.  The
subspace $C^\infty(H_\pi) \subset H_\pi$ is endowed with the
Fr\'echet $C^\infty$ topology, that is, the topology defined
by the family of seminorms $\{ \|\cdot \|_{E_1,E_2,
  \dots,E_m} \,\vert \, m \in \N \text{ and
}E_1,\dots,E_m \in \nil\}$ defined as follows:
$$
\|v\|_{E_1,E_2, \dots,E_m} := \|\pi_*(E_{1})
\pi_*(E_{2})\cdots \pi_*(E_{m}) v\| , 
\quad \text{ for }v\in C^\infty(H_\pi)\,.
$$
\begin{definition}
  Let $ \pi $ be a unitary representation of the Lie algebra $\nil$ on
  the Hilbert space $H_{\pi} $.  The space of\/ {\em Schwartz
    distributions for the representation $ \pi $} is defined as the
  the dual space 
  of $ C^\infty(H_ \pi)$ (endowed with the Fr\'echet $C^\infty$
  topology). The space of Schwartz distributions for the
  representation $ \pi $ will be denoted ${\Cal D}' (H_{\pi})$.
\end{definition}

The representation $ \pi_*$ extends in a canonical way to a
representation (denoted by the same symbol) of the
enveloping algebra $\Cal U (\nil)$ of $\nil$ on the Hilbert
space $H_\pi$. Let $\Delta\in \Cal U(\nil)$ a left-invariant
second-order, positive elliptic operator on $\Nil$ fixed
once for all. For example, the operator $\Delta= -(V_1^2 +
\dots + V_d^2)$, where $\{V_1,\dots, V_d\}$ is a basis
of~$\nil$ as a vector space.  For all $ \alpha\in \R^{+} $,
let $W^{\alpha}(H_{\pi}) \subset H_{\pi} $ be the Sobolev
space of vectors in the maximal domain of the essentially
self-adjoint operator $ \left(I+
  \pi_*(\Delta)\right)^{\alpha/2} $ endowed with the Hilbert
space norm
$$
\| v \|_\alpha\, := \, \| \left(I+
  \pi_*(\Delta)\right)^{\alpha/2} v \|\,, \quad \text{for
  all } v\in W^{\alpha}(H_{\pi})\,.
$$
Let also $ W^{-\alpha}(H_{\pi}) \subset {\Cal D}' (H_{\pi})
$ be the dual Hilbert space of $ W^{\alpha}(H_{\pi}) $. Then
$C^\infty(H_\pi)$ is the projective limit of the spaces $
W^{\alpha}(H_{\pi})$ (and consequently ${\Cal D}' (H_{\pi})$
is the inductive limit of $W^{-\alpha}(H_{\pi})$) as
$\alpha\to +\infty $. A distribution $D \in W^{-\alpha}(H_{\pi}) $ 
will be called a {\em distribution of order at most} $\alpha\in \R^+$.

\begin{definition}  
  Let $\pi $ be a unitary representation of the Lie group
  $\Nil $ (with Lie algebra $\nil$) on the Hilbert space $
  H_{\pi} $. The \/{\em cohomological equation for $ X\in
    \nil $} in the representation $ \pi $ is the linear
  equation
  \begin{equation}
    \pi_*(X)u=f \quad (\text{for the unknown vector } u\in H_\pi),
  \end{equation}  
  where $ f\in H_{\pi} $ is a given (sufficiently smooth)
  vector.
\end{definition}

\begin{definition}  
  For any $X\in \nil$, the space of \/{\em $ X $-invariant
    distributions }for the representation $\pi$ is defined
  as the space $ {\Cal I}_{X}(H_{\pi } ) $ of all
  distributional solutions $ D \in {\Cal D}' (H_{\pi})
  $ of the equation $ \pi_*(X) D=0 $.  Let
  $$ 
  {\Cal I}^{\alpha}_ {X}(H_{\pi } ):={\Cal I}_{X}(H_{\pi }
  )\cap W^{-\alpha}(H_{\pi})
  $$ 
  be the subspace of invariant distributions of order at
  most $ \alpha\in \R^{+}$.
\end{definition}

It is immediate that all $ X $-invariant distributions in
$\Cal D'(H_\pi)$ yield obstructions to the existence of
$C^{\infty}$ solutions $u\in C^\infty(H_{\pi})$ of the
cohomological equation and that $X$-invariant distributions
of order at most $\alpha\in\R^+$ yield obstructions to the
existence of solutions $u\in W^{\alpha+1}(H_\pi)$.

Let $\Cal O$ be a codjoint orbit of maximal rank and let
$\pi \in \Pi_{\Cal O}$ be an induced irreducible (adapted)
unitary representation of $\Nil$ on $H_\pi= L^2(\R, H')$
constructed as in Lemma~\ref{lemma:cohomeqone}. The space of
$C^\infty(H_\pi)$ of $C^\infty$ vectors for the
representation $\pi$ can be characterized as follows.

Let $\Cal S\left(\R, C^\infty(H')\right)$ be the space of
Schwartz functions on $\R$ with values in the Fr\'echet
space of $C^\infty$ vectors for the representation $\pi'$ on
$H'$ (cf.~Lemma~\ref{lemma:cohomeqone}). The space $\Cal
S\left(\R, C^\infty(H')\right)$ is by definition endowed with the
Fr\'echet topology induced by the family of seminorms
$$\{  \Vert \cdot \Vert _{i,j,E_1,E_2,\dots, E_m} \,\vert\, i,j,m\in \N
\text{ and }E_1, \dots,E_m\in\nil'\}\, $$ 
defined as follows: for all $f\in \Cal S \left(\R,
  C^\infty(H')\right)$,
\begin{equation}
\label{eq:seminorms} 
\Vert f \Vert _{i,j,E_1,E_2,\dots, E_m} := \sup_{t\in \R}
\Vert (1+ t^2)^{j/2}\,\pi'_*(E_1) \cdots \pi'_*(E_m)
f^{(i)}(t)\Vert_{H'} \,.
\end{equation}

\begin{lemma}
\label{lemma:smoothvectors} 
As topological vector spaces,
\begin{equation*}
C^\infty(H_\pi) \equiv  \Cal S\left(\R, C^\infty(H')\right)\,.
\end{equation*}
\end{lemma}
\begin{proof}
  Let $X$, $Y\in \nil$ and $\nil'\subset \nil$ be as in
  Lemma~\ref{lemma:cohomeqone}. By formula~\pref{eq:semi3},
  for any $E\in \nil'$ and for all $f\in \Cal S\left(\R,
    C^\infty(H')\right)$,
\begin{equation}
\label{eq:semi4}
\left(\pi_*(E)f \right)(t)= 
\sum_{j=0}^{k} \frac{t^j}{j!} \pi'_* 
\left( \text{\rm ad}(X)^j E\right)f(t).
\end{equation}
We claim that the Schwartz space $\Cal S\left(\R,
  C^\infty(H')\right) \subset C^\infty(H_\pi)$ and that the
embedding is continuous.  In fact, for any $f\in \Cal
S\left(\R, C^\infty(H')\right)$, we have
\begin{equation}
\label{eq:Sboundzero}
\Vert f \Vert \le
\sup_{t\in \R} \Vert (1+ t^2)^{1/2}\, f (t)\Vert_{H'}\,,
\end{equation}
hence it follows from formula~\pref{eq:semi4} that there
exists a constant $C_k>0$ such that for any $E\in \nil'$ and
all $f\in \Cal S\left(\R, C^\infty(H')\right)$,
\begin{equation}
\label{eq:Sboundone}
\Vert \pi_*(E) f \Vert  \leq  C_k  \max_{0\leq j\leq k} 
\left[ \sup_{t\in \R}   \Vert (1+t^2)^{\frac{1}{2}} t^j
\pi'_* \left( \text{\rm ad}(X)^j E\right)f(t)\Vert_{H'} \right] \,.
\end{equation}
In order to simplify the notation, let $E^{(j)}:= \text{\rm
  ad}(X)^j E \in \nil'$ for any $E\in \nil'$ and for any
$j\in\{0, \dots, k\}$. By estimate~\pref{eq:Sboundone} and
by iterated applications of the identity~\pref{eq:semi4}, it
follows that for any $\ell\in \N$ there exists a constant
$C_{k,\ell}>0$ such that, for any $E_1, \dots, E_\ell \in
\nil'$ and for all $f\in \Cal S\left(\R,
  C^\infty(H')\right)$,
\begin{multline}
\label{eq:Sboundtwo}
\Vert \pi_*(E_1)\cdots \pi_*(E_\ell) f \Vert \leq \\ \leq
C_{k,\ell} \, \max_{0\leq j_1, \dots, j_\ell \leq k}\,
\left[ \sup_{t\in \R} \Vert (1+ t^2) ^{\frac{1}{2}}
  \prod_{\alpha=1}^{\ell} t^{j_\alpha} \pi'_* \left(
    E^{(j_\alpha)}_\alpha\right) f(t)\Vert_{H'} \right],
\end{multline}
hence, for any $E_1, \dots, E_\ell \in \nil'$, any $i\in\N$
and any $f\in \Cal S\left(\R, C^\infty(H')\right)$,
\begin{multline}
\label{eq:Sboundthree}
\Vert \pi_*(E_1)\cdots \pi_*(E_\ell) \pi_*(X)^i f \Vert \leq
\\ \leq C_{k,\ell} \, \max_{0\leq j_1, \dots, j_\ell \leq
  k}\, \left[ \sup_{t\in \R} \Vert (1+ t^2) ^{\frac{1}{2}}
  \prod_{\alpha=1}^\ell t^{j_\alpha} \pi'_* \left(
    E^{(j_\alpha)}_\alpha\right) f^{(i)}(t)\Vert_{H'}
\right].
\end{multline}
Since $[X, \nil']\subset \nil'$ the claim follows.

\smallskip Conversely, we claim that $C^\infty(H_\pi)
\subset \Cal S\left(\R, C^\infty(H')\right) $.  In fact, let
$f\in C^\infty( H_\pi)$ be a $C^\infty$ vector. Since
$\pi_*(X)$ is the operator $f\mapsto f'$, we have $f\in
C^\infty(\R, H')$. Since $\pi_*(Y)$ equals, up to a non-zero
(scalar) factor, the pointwise operator $f(t) \mapsto t
f(t)$, we have $(t\mapsto t^j f(t))\in L^2(\R, H')$, for all
$j\ge 0$, hence $f\in \mathcal S(\R, H')$. Thus
\begin{equation}
\label{eq:inclusion1}
C^\infty(H_\pi) \subset \mathcal S(\R, H')
\end{equation} 
and we therefore have a continuous linear operator 
$$\mathcal E_t: C^\infty(H_\pi) \to H'\,, \qquad \mathcal E_tf = f(t)\,.$$

For any $E\in \nil'$ and any $j \in \N$, let $T^j\otimes
\pi'_*(E)$ be the densely defined linear operator on
$L^2(\R, H')$ uniquely determined by the identities
 $$
 {\Cal E}_t [T^j\otimes \pi'_*(E)] f := t^j\pi'_*(E)
 {\mathcal E}_t f\,, \quad \text{ for all }t\in \R\,.
 $$
 We will prove by (finite) induction that, for all
 $\alpha\in \{0,\dots, k-1\}$, if $E\in \nil_{k-\alpha}\cap
 \nil'$, then we have, for all $j \in \N$ and for all $f\in
 C^\infty(H_\pi)$,
\begin{equation}
  \label{eq:piprime1}
  \begin{aligned}
 & f \in \dom[T^j\otimes \pi'_*(E) ]  \, ,\\
   [&T^j \otimes \pi'_*(E)] f \in  \mathcal C^{\infty}(H_\pi)\,.
    \end{aligned}
    \end{equation}

    Remark that, since $[T^j\otimes \pi'_*(E)] $ equals
    $\pi_*(Y)^j [1\otimes\pi'_*(E)] $ up to a non-zero
    scalar factor. the second line of formula
    $\pref{eq:piprime1}$ actually follows from
\begin{equation}
  \label{eq:piprime2}
   [1\otimes \pi'_*(E)] f \in  \mathcal C^{\infty}(H_\pi)\,.
\end{equation}

If $\alpha=0$, for any $E_0 \in \nil_k \subset Z(\nil)$ the
operators $\pi_*(E_0)$ and $\pi'_*(E_0)$ are scalar
multiples of the identity. Hence the claim
\pref{eq:piprime1} holds in this case.

Let us assume \pref{eq:piprime1} for all $E \in
\nil_{k-\alpha}$ and for all $j\in \N$. Let $E_{\alpha+1}
\in \nil_{k-(\alpha+1)}\cap \nil'$. Remark that the elements
$$
E^{(\beta)}_\alpha:=\text{\rm ad}(X)^{\beta} E_{\alpha+1}
\in \nil_{k-\alpha}\cap \nil' \,, \quad \text{ for all }\,
\beta\in \{1,\dots, \alpha+1\}\,.
$$ 
By formula~\pref{eq:semi4} we have the identity,
\begin{equation}
\label{eq:piprimeid} 
\pi'_*(E_{\alpha+1}){\Cal E}_t f =  
{\Cal E}_t \pi_*(E_{\alpha+1}) f  -
\sum_{\beta=1}^{\alpha+1} \frac{1}{\beta!} 
{\Cal E}_t [T^\beta\otimes \pi' _*( E^{(\beta)}_\alpha)] f\,.
\end{equation}
where $\pi_*(E_{\alpha+1}) f \in C^{\infty}(H_\pi)$ and
$[T^i\otimes \pi' _*( E^{(\beta)}_\alpha)] f \in
C^\infty(H_\pi)$ by the induction hypothesis, since $f\in
C^\infty(H_\pi)$. It follows by the inclusion
\pref{eq:inclusion1} and the identity \pref{eq:piprimeid}
that the claim \pref{eq:piprime1} holds for all
$E_{\alpha+1} \in \nil_{k-\alpha}\cap \nil'$ and all $j\in
\N$. The induction step is therefore completed.

It follows from \pref{eq:piprime1} by another induction
argument (on $i\in \N$), that, for all $E\in \nil'$ and all
$i$, $j \in N$, we have
\begin{equation}
  \label{eq:piprime3}
  \begin{aligned}
 & f^{(i)} \in \dom[T^j\otimes \pi'_*(E) ]   \, ,\\
   [&T^j \otimes \pi'_*(E)] f^{(i)}  \in  \mathcal C^{\infty}(H_\pi)\,.
    \end{aligned}
\end{equation}
In fact, for $i=0$ the conditions \pref{eq:piprime3} reduce
to \pref{eq:piprime1}. Assume that \pref{eq:piprime3} holds
for $i\in \N$ and for all $j\in \N$.  Since $[1\otimes
\pi'_*(E) ] f^{(i)} \in C^{\infty}(H_\pi)$, we have
$$
[1\otimes \pi'_*(E)] f^{(i+1)} =  
\pi_*(X) [1\otimes \pi'_*(E) ] f^{(i)}  
\in C^{\infty}(H_\pi) \subset \Cal S(\R,
H')\,,
$$
hence $ f^{(i+1)}\in \dom[1\otimes \pi'_*(E)] $. Next, for
$j \in \N\setminus \{0\}$, we have the following immediate
identity:
\begin{equation}
[T^j \otimes \pi'_*(E)] f^{(i+1)} = 
\pi_*(X) [T^j \otimes \pi'_*(E)] f^{(i)} -
 j [T^{j-1} \otimes \pi'_*(E)] f^{(i)} \,,
\end{equation}
which allows to conclude that \pref{eq:piprime3} holds for
$i+1$. The induction argument is therefore completed.

Finally, by successive applications of~\pref{eq:piprime3} we
have that, for any $i$, $j\in \N$, any $E_1$, \dots, $E_m\in
\nil'$ and for all $f\in C^\infty(H_\pi)$,
\begin{equation}
  \label{eq:piprime4}
  \begin{aligned}
    & f^{(i)} \in \dom[ T^j \otimes \pi'_*(E_1) \cdots \pi'_*(E_m)]\,,  \\
    [&T^j \otimes \pi'_*(E_1) \cdots \pi'_*(E_m)] f^{(i)}
    \in \mathcal C^{\infty}(H_\pi)\,.
    \end{aligned}
\end{equation}
 
\smallskip In conclusion, we have proved that $\Cal
S\left(\R, C^\infty(H')\right) \subset C^\infty(H_\pi)$ and
that the inclusion is continuous and surjective, hence it is
an isomorphism of Fr\'echet spaces by the open mapping
theorem (or by a direct argument based on the proof of
surjectivity given above).
\end{proof}

\subsection{Main Estimates} 
In this section we describe all invariant distributions and
we prove Sobolev bounds for the Green operator of the
cohomological equation in each irreducible unitary
representation of maximal rank (determined by a coadjoint
orbit $\Cal O$ of maximal rank) for any admissible vector
$X\in \nil\setminus \nil_{k-1}^\perp(\Cal O)$.

Let $\< \cdot,\cdot\>$ be an euclidean product on the Lie
algebra $\nil$ fixed once for all and let $\Vert \cdot
\Vert$ denote the corresponding norm.

Let~$ \Cal O $ be any coadjoint orbit.  Since the
restrictions $\lambda\vert\nil_k$ and $ B_{\lambda} \vert
\nil\times \nil_{k-1} $ do not depend on the choice of the
linear form $ \lambda\in \Cal O $, the following definitions
are well-posed.

\smallskip Let
\begin{equation}
  \label{eq:ws}
  \begin{aligned}
    w_k(\Cal O):=& \Vert \lambda\vert \nil_k\Vert =
    \max_{Z\in \nil_k\setminus\{0\}} \frac{\vert
      \lambda(Z)\vert}{\Vert Z\Vert} \,;\\
    w_Z(\Cal O):=& \Vert \lambda\vert Z(\nil)\Vert =
    \max_{Z\in Z(\nil)\setminus\{0\}} \frac{\vert
      \lambda(Z)\vert}{\Vert Z\Vert} \,.
  \end{aligned}
\end{equation}
Since $\nil_k\subset Z(\nil)$ we have $w_k(\Cal O) \leq
w_Z(\Cal O)$.

\smallskip For all $(X, Y) \in \nil\times \nil_{k-1} $, let
\begin{equation}
  \label{eq:delta}
  \begin{aligned}
    &\delta_{\Cal O}(X,Y):=\vert B_{\lambda}(X, Y) \vert
    \,\,,
    \,\,\text{for any } \lambda\in \Cal O \,, \text{ and} \\
    &\delta_{\Cal O}(X):=\max \{ \delta_{\Cal O}(X,Y) \,\mid
    \,Y \in \nil_{k-1} \text{ and } \Vert Y \Vert =1 \} \,.
  \end{aligned}
\end{equation}

By definition we have $ \delta_{\Cal O}(X)>0 $ if and only
if~$ X \not\in \nil_{k-1} ^{\perp}(\Cal O) $.  The latter
condition is non-empty if and only if $ \Cal O $ has maximal
rank, in which case it holds except for a subspace of
positive codimension.

\begin{lemma}
  \label{lem:invdist1} 
  Let $ \delta_{\Cal O}(X)>0 $, and $ \pi\in \Pi_{\Cal O} $
  an irreducible essentially skew-adjoint representation of
  $ \nil $ on a Hilbert space $ H_{\pi} $. Let ${\Cal I}_X(
  H_{\pi})\subset {\Cal D}'(H_{\pi})$ be the space of $ X
  $-invariant distributions, i.e. of linear functionals
  which vanish on $ \text{\rm Ran}\left(X\vert
    C^\infty(H_{\pi})\right)\subset C^\infty(H_{\pi}) $.
  The space ${\Cal I}_X(H_{\pi})$ has a countable basis
  $\Cal B_X(H_\pi)$ with the following properties:
  \begin{enumerate}
  \item each $D\in {\Cal B}_{X}(H_{\pi})$ has Sobolev order
    equal to $1/2$, in the sense that $D \in
    W^{-\alpha}(H_\pi)$ for any $\alpha>1/2$;
  \item for each $D\in {\Cal B}_{X}(H_{\pi})$ and each $
    \alpha>1/2 $, there exists a constant $ C:=C(D,
    \alpha)>0 $ such that the following holds: for any $
    Y\in \nil_{k-1} $ with $ \delta_{\Cal O}(X,Y) >0 $ and
    for all $ f\in W^{\alpha}(H_ {\pi}) $, we have
  \begin{equation}
    \label{eq:invdistnorm} 
    \vert D(f)\vert \le C\,
    \delta_{\Cal O}(X,Y) ^{-1/2}\, \Vert (I- \pi_*(Y)^{2}) ^{\alpha/2}f
    \Vert \,\,.
\end{equation}
\end{enumerate}
\end{lemma}


\begin{proof}
  Let $ \Cal O\subset \nil^{\ast} $ be a coadjoint orbit
  such that $ \delta_{\Cal O}(X)>0 $, let $ Y\in \nil_{k-1}
  $ be such that $ \delta_{\Cal O}(X,Y)>0 $ and let $ \pi
  \in \Pi_{\Cal O} $.  Since the statement of the theorem
  depends only on the unitary equivalence class of the
  representation, we can assume by Lemma~
  \ref{lemma:cohomeqone} that~$ \pi $ is obtained from an
  irreducible representation~$ \pi' $ of an ideal $
  \nil'\subset \nil $ of codimension~$ 1 $ on a Hilbert
  space $ H' $ by adjoining~$ X\in \nil\setminus \nil' $,
  hence $ \pi_*(X) $ is the derivative operator $
  \frac{\partial}{ \partial t } $ on $ H_{\pi}\equiv
  L^{2}(\R,H', dt) $, and that $ \pi_*(Y) = 2\pi\imath \,t
  \,\delta_{\Cal O}(X,Y)\, \I_{H'} $.

  Let  $\Cal B' \subset H'$ be a basis of the Hilbert space $H'$
  and let ${\Cal B}_{X}(H_{\pi})$ be the set of linear
  functionals $D_e$, $ e\in \Cal B'$, defined as
  follows. For any $e\in B'$, for all  $f\in  \mathcal S(\R,H')$, let
  \begin{equation}
    \label{eq:invdist} D_{e}(f) := \int_{\R} \<f(t),e\>_
    {H'} \,dt \,\, .
  \end{equation}
  
  By H\"older inequality, for all $ \alpha>1/2 $ there
  exists a constant $ C_{\alpha}>0 $ such that the following
  (a priori) estimate holds:
  \begin{equation}
    \label{eq:invdistest} \vert D_{e}(f) \vert \leq
    \frac{C_{\alpha}}{\delta_{\Cal O}(X,Y)^{1/2}}\, \vert e\vert_
    {H'} \, \Vert \left(I- \pi_*(Y)^{2}\right)^{\alpha/2} f
    \Vert \,\,;
  \end{equation}
  hence $ D_{e} \in W^{-\alpha}(H) $ for all $
  \alpha>1/2 $; in particular $D_{e}\in \Cal D'(H_\pi)$.

  Clearly $D_{e}$ is invariant by translations, hence
  $D_e \in {\Cal I}_{X}(H_{\pi})$.  It remains to be
  proved only that ${\Cal B}_{X}(H_{\pi})$ is a basis of
  ${\Cal I}_{X}(H_{\pi})$.  Since by
  Lemma~\ref{lemma:smoothvectors} we have
  $C^{\infty}(H_\pi) \equiv \Cal S\left(\R,
    C^{\infty}(H')\right)$, if $D_e(f)=0$ for all $e\in
  \Cal B'$, then $\int_\R f(s) ds =0 \in C^{\infty} (H')$
  and the function
  $$
  G(t):= \int_{-\infty}^t f(s)\, ds \,\, \in\,\, \Cal
  S\left(\R, C^{\infty}(H')\right)\,.
  $$
  It follows that $D (f)=0$ for all $D\in {\Cal
    I}_{X}(H_{\pi})$, since there exists $g\in
  C^\infty(H_\pi)$ such that $f= \pi(X)g$. Hence, by the
  Hahn-Banach theorem, ${\Cal B}_{X}(H_{\pi})$ is a
  (countable) system of generators for the (closed) subspace
  ${\Cal I}_{X}(H_{\pi}) \subset \Cal D'(H_\pi)$. It is
  immediate to verify that, since $\Cal B' \subset H'$ is
  linearly independent, the system ${\Cal B}_{X}(H_{\pi})$
  is also linearly independent, hence it is a basis of
  ${\Cal I}_{X}(H_{\pi})$.
\end{proof}
\medskip
\begin{theorem}
  \label{th:cohomeq} Let $ \delta_{\Cal O}(X)>0 $, and let
  $\pi\in \Pi_{\Cal O} $ be an irreducible essentially
  skew-adjoint representation of $ \nil $ on a Hilbert space
  $ H_{\pi} $. 
  \begin{enumerate}
  \item There exists a Green operator $ G_{X}:
    W^{\alpha}(H_{\pi})\to W^{\beta}(H_{\pi}) $, defined for
    all $ \alpha>1/2 $ and $ \beta<-1/2 $, such that, for
    all $ f\in W^{\alpha} (H_{\pi}) $, the distribution $
    u:=G_{X}f $ is a solution of the cohomological equation
    $ \pi_*(X)u=f $; there exists a constant $
    C:=C(X,\alpha)>0 $ such that
    \begin{equation*}
      \Vert G_{X}f \Vert _{\beta} \le C\, \delta_{\Cal O}(X)^
      {-1}\, \Vert f \Vert _{\alpha} \,.
    \end{equation*}
  \item If $ f\in W^{\alpha}(H_{\pi}) $, $ \alpha>1 $, and $
    D(f)=0 $, for all $ D\in {\Cal
      I}_{X} (H_{\pi}) $, then $ G_{X}f \in
    W^{\beta}(H_{\pi}) $, for all $ \beta<(\alpha-1)/k $;
    furthermore there exists a constant $
    C':=C(X,\alpha,\beta)>0 $ such that
    \begin{equation*}
      \Vert G_{X}f \Vert _{\beta} \le C'\, 
      \max\{ 1,w_k(\Cal O)^\beta\}\, \max \{1,\delta_
      {\Cal O}(X)^{-1-k \beta} \}\, \Vert f \Vert _{\alpha}
      \,\,.
    \end{equation*}
  \end{enumerate}
\end{theorem}
\begin{proof}
  Let $ \Cal O\subset \nil^{\ast} $ be a coadjoint orbit and
  let $ \pi \in \Pi_{\Cal O} $.  If $ \delta_{\Cal
    O}:=\delta_ {\Cal O}(X)>0 $, by
  Lemma~\ref{lemma:cohomeqone} we can assume that the
  representation~$ \pi $ is obtained from an irreducible
  representation~$ \pi' $ of a codimension~$ 1 $ ideal $
  \nil' $ on a Hilbert space $ H' $ by adjoining~$ X\in
  \nil\setminus \nil' $, hence $ \pi_*(X) $ is the
  derivative operator $ \frac{\partial} { \partial t } $ on
  $ H_{\pi}\equiv L^{2}(\R,H', dt) $, and that there exists
  $ Y\in \nil_{k-1} $ such that $ \pi_*(Y) = 2\pi\imath \,t
  \,\delta_{\Cal O}\, \I_{H'} $.  In fact, by compactness of
  the unit sphere in $ \nil_{k-1} $, there exists $ Y\in
  \nil_{k-1} $ such that $ \Vert Y \Vert=1 $ and $
  \delta_{\Cal O}(X,Y)= \delta_{\Cal O}>0 $.

  \smallskip
  \noindent
  The operator $ G_{X}: C_{0}(\R,H') \to B(\R,H') $ defined
  for $ f\in C_{0}(\R,H') $ by
  \begin{equation}
    \label{eq:Greenopt} G_{X}f (t) := \int_{-\infty}^{t} f(s)\,
    ds \,\,, \quad \text{for all } t\in \R\,,
  \end{equation}
  admits a bounded extension $ G_{X}: W^{\alpha}(H_{\pi})
  \to C_{B}(\R,H') $, for all $ \alpha>1/2 $; in fact , if $
  f\in W^{\alpha}(H_{\pi}) $, $ \alpha>1/2 $, by H\"older
  inequality there exists a constant $ C_{\alpha}>0 $ such
  that
  \begin{equation}
    \int_{-\infty}^{t} \vert f(s)\vert_{H'}\, ds \,\, \leq
    \frac{C_{\alpha}}{\delta_{\Cal O}^{1/2} }\, \Vert \left(I-\pi
      (Y)^{2}\right)^{\alpha/2} f \Vert\,.
  \end{equation}
  A similar estimate shows that $ C_{B}(\R,H')\subset
  W^{-\alpha} (H_{\pi}) $ for all $ \alpha>1/2 $.  In fact,
  if $ u\in C_ {B}(\R,H') $, for all $ \alpha>1/2 $ there
  exists $ C_{\alpha}>0 $ such that, if $ v\in
  W^{\alpha}(H_{\pi}) $,

  \begin{equation}
    \vert (u,v) \vert \leq \int_{\R} \vert\<u(s),v(s)\>_{H'}
    \vert \, ds \,\leq \, \frac{C_{\alpha}}{\delta_{\Cal O}^{1/2} }\,
    \vert u\vert_{\infty} \, \Vert v \Vert_{\alpha} \,\,.
  \end{equation}
  It follows that $ G_{X}: W^{\alpha}(H_{\pi}) \to W^{\beta}
  (H_{\pi}) $ is well-defined and bounded for all $
  \alpha>1/2 $ and all $ \beta<-1/2 $, and there exists $
  C_{\alpha,\beta}>0 $ such that, for all $ f\in
  W^{\alpha}(H_{\pi}) $,

  \begin{equation}
    \label{eq:Greenestone} \Vert G_{X}f \Vert _{\beta}\,\,
    \leq \,\, \frac{C_{\alpha,\beta}}{\delta_{\Cal O}} \,\,
    \Vert f \Vert_{\alpha}
  \end{equation}
  By construction $ \, \pi_*(X) G_{X}f = f\, $ in $
  W^{\beta-1} (H_{\pi}) $, for all $ f\in
  W^{\alpha}(H_{\pi}) $, $ \alpha>1/2 $ and $ \beta<-1/2 $,
  and $ \,G_{X} \pi_*(X)u=u\, $ in $ W^{\beta}(H_{\pi}) $,
  for all $ u\in W^{\alpha+1}(H_{\pi}) $, $ \alpha>1/2 $ and
  $ \beta<-1/2 $.

  \smallskip Finally, if $ f\in W^{\alpha}(H_{\pi}) $, for $
  \alpha>1/2 $, and $ D(f)=0 $, for all $ D \in
  {\Cal I }^{\alpha}_ {X}(H_{\pi}) $,

  \begin{equation}
    G_{X}f (t) := \int_{-\infty}^{t} f(s)\, ds \,\,= \,\, -
    \int_{t}^{+\infty} f(s)\, ds \,\, , \quad \text{for all
    } t\in \R\,.
  \end{equation}

  It follows that

  \begin{multline}
    \Vert \pi_*(Y)^{\ell} G_{X}f \Vert ^{2} \leq
    \int_{0}^{+\infty} \left( \vert 2\pi \delta_{\Cal
        O}\,t\vert^{\ell}\, \int_ {t}^{+\infty} \vert
      f(s)\vert_{H'}\, ds \right)^{2}\,dt
    \\
    + \,\, \int_{-\infty}^{0} \left(\vert 2\pi \delta_{\Cal
        O}\,t\vert^{\ell}\, \int_{-\infty}^{t} \vert
      f(s)\vert_{H'}\, ds \right)^{2}\,dt\, \,\,.
  \end{multline}
  Since, for all $ \alpha>1 $,
  \begin{multline}
    C^{2}_{\alpha,\ell} := \int_{0}^{+\infty} (2\pi
    t)^{2\ell} \left( \int_{t}^{+\infty} ( 1 + 4\pi
      s^{2})^{-(\ell+\alpha)}\,
      ds \right)dt \\
    \,\,+ \,\, \int_{-\infty}^{0} (2\pi
    t)^{2\ell}\left(\int_ {-\infty}^{t} ( 1 + 4\pi
      s^{2})^{-(\ell+\alpha)}\, ds \right)dt \,\,< \,\,
    +\infty \,\,,
  \end{multline}
  by H\"older inequality and change of variables,
  \begin{equation}
    \label{eq:Greenesttwo} \Vert  \pi_*(Y)^{\ell} G_{X}f \Vert
    \leq \frac{C_{\alpha,\ell}}{\delta_{\Cal O}}\, \Vert
    \left(I- \pi_*(Y)^{2}\right)^{\frac{\ell+\alpha}{2}} f\Vert
    \,.
  \end{equation}
  
  Let $ {\mathfrak U }(\nil) $ and $ {\mathfrak U }(\nil') $
  be the enveloping algebras of~$ \nil $ and~$ \nil' $ respectively.
  Clearly from the inclusion $\nil'\subset\nil$ we have that $ {\mathfrak
    U }(\nil') $ is contained in $ {\mathfrak U }(\nil) $ as
  a subalgebra.  In addition, since $\nil'$ is an ideal
  of~$\nil$, $(\Ad (\exp t X))_{t\in \R}$ is a one-parameter
  group of automorphisms of $\nil'$ which extends to a 
  one-parameter group $(\exp t \delta)_{t\in \R}$ of automorphisms 
  of $ {\mathfrak U }(\nil') $. The generator $\delta$ of this group of
  automorphisms is the derivation on $
  {\mathfrak U }(\nil') $ obtained by extending the
  derivation~$ \ad (X) $ of~$ \nil' $ from $ \nil' $ to $
  {\mathfrak U }(\nil') $. Observe that from the nilpotency
  of $\nil$ it follows that, for any $P\in {\mathfrak U
  }(\nil')$ there exists a first integer $[P]$ such that
  $\delta^{[P]+1} P=0$.
  
  \smallskip We recall that by
  Lemma~\ref{lemma:smoothvectors} we have
  $C^\infty(H_\pi)=\mathcal S(\R,C^\infty(H'))$. Let
  ${\mathcal E}_t: C^\infty(H_\pi) \to C^\infty(H')$ be the
  linear (evaluation) operator defined by ${\mathcal E}_t f
  = f(t)$. By definition of the representation $\pi_*$ on
  $C^\infty(H_\pi)$ (cf. formula~\pref{eq:semi3}), we have,
  for any $P\in {\mathfrak U }(\nil') $ and any $E\in
  \nil'$,
  \begin{equation}
    \label{eq:piaction}
    {\mathcal E}_t \, \pi_*(P) =  
    \pi_*'(e^{t \delta} P) \, {\mathcal E}_t= \sum_{j=0}^{[P]}
    \frac {t^j} {j!}  \pi_*'(\delta^j P)
    \, {\mathcal E}_t, \quad \text{for all } \,t\in \R\,,
  \end{equation}
  from which we obtain, for all $s$, $t\in \R$,
  \begin{equation}
    \label{eq:piprimaction}
    \pi_*'(e^{(t-s )\delta} P) \, {\mathcal E}_t=
    {\mathcal E}_t \, \pi_*(e^{-s \delta} P) = 
    \sum_{j=0}^{[P]}
    \frac {(-s)^j} {j!} {\mathcal E}_t \,  
    \pi_*(\delta^j P)\,. 
  \end{equation}
  Since the action of $\pi_*(Y)$ on $C^\infty(H_\pi)$ can be rewritten as
  \begin{equation}
    \label{eq:pY}
    {\mathcal E}_t \pi_*(Y^j) = \left( 2\pi \imath \delta_{\mathcal
        O} t\right)^j {\mathcal E}_t \,, 
  \quad \text{for all } \,j\in\N \text{ and } t\in \R\,,
  \end{equation}
  setting $t=s$ in the \pref{eq:piprimaction} we obtain, for
  all $t\in \R$,
  \begin{equation}
    \label{eq:piprimactionbis}
    \pi_*'( P) \, {\mathcal E}_t=
    {\mathcal E}_t \,\pi(e^{-t \delta} P) = \sum_{j=0}^{[P]}
    \frac {(-1)^j} {j!(2\pi \imath \delta_{\mathcal O})^j} 
    {\mathcal E}_t \, \pi_*(Y^j \delta^j P) \,.
  \end{equation}
  
  Since ${\mathcal E}_t G_X=\int_{-\infty}^t {\mathcal E}_s
  \,ds$ we have, for all $t\in\R$,
  
 \begin{equation*}
 \begin{aligned}
    &{\mathcal E}_t \pi_*(P) G_X = \pi_*'(e^{t\delta}P) {\mathcal E}_t G_X=
    \pi_*'(e^{t\delta}P) \int_{-\infty}^t {\mathcal E}_s
    \,ds \\ 
    &=\sum_{j=0}^{[P]} \frac {t^j} {j!}
    \pi_*'(\delta^j P) \int_{-\infty}^t {\mathcal E}_s \,ds=
    \sum_{j=0}^{[P]} \frac {t^j} {j!}  \int_{-\infty}^t
    \pi_*'(\delta^j P) {\mathcal E}_s
    \,ds \\
   &= \sum_{j=0}^{[P]} \frac {t^j} {j!}  \sum_{\ell=0}^{[P]}
    \frac {(-1)^\ell} {\ell!(2\pi \imath \delta_{\mathcal
        O})^\ell}\int_{-\infty}^t {\mathcal E}_s \, \pi_*(Y^\ell
    \delta^{j+\ell} P)
    \,ds\\
    &=\sum_{j=0}^{[P]}\sum_{\ell=0}^{[P]} \frac {t^j} {j!}
    \frac {(-1)^\ell} {\ell!(2\pi \imath \delta_{\mathcal
        O})^\ell} {\mathcal E}_t
    G_X \pi_*(Y^\ell \delta^{j+\ell} P)\\
    &={\mathcal E}_t \sum_{j=0}^{[P]}\sum_{\ell=0}^{[P]} \frac
    {(-1)^\ell} {j!\ell!(2\pi \imath \delta_{\mathcal
        O})^{j+\ell}} \pi_*(Y^j)
    G_X \pi_*(Y^\ell \delta^{j+\ell} P)\\
    &={\mathcal E}_t \sum_{m=0}^{[P]}\frac {1} {(2\pi \imath
      \delta_{\mathcal O})^m}\sum_{\ell+j =m} \frac
    {(-1)^\ell} {j!\ell!} \pi_*(Y^j) G_X \pi_*(Y^\ell
    \delta^{m} P)\,,
    \end{aligned}
    \end{equation*}
hence, using~\pref{eq:Greenesttwo}, for all $ \alpha>1 $,
  there exists a constant~$ C_{\alpha,[P]}>0 $ such that for
  all $f\in C^{\infty}(H_\pi)$, we have
  \begin{multline}
    \label{eq:bound} \Vert \pi_*(P) G_{X}f \Vert\le \\
    \sum_{m=0}^{[P]}\frac {1} {(2\pi \delta_{\mathcal
        O})^m}\sum_{\ell+j =m} \frac {1}
    {j!\ell!}\frac{C_{\alpha,j}}{\delta_{\Cal O}}\, \Vert
    \left(I-\pi(Y)^{2}\right)^{\frac{j+\alpha}{2}}
    \pi_*(Y^\ell
    \delta^{m} P) f\Vert\le \\
    C_{\alpha,[P]} \sum_{m=0}^{[P]} {\delta_{\mathcal
        O}^{-m-1}} \Vert
    \left(I-\pi(Y)^{2}\right)^{\frac{m+\alpha}{2}} \pi_*(
    \delta^{m} P) f\Vert.
  \end{multline}
  
  Let $\Delta$ the positive definite Laplacian associated to
  an euclidean product $\<\cdot,\cdot\>$ on the Lie
  algebra~$\nil$. Then for any orthonormal basis $\mathcal
  B$ of~$\nil$ we have $\Delta=- \sum_{V\in \mathcal B}
  V^2$.  Let $U\in\nil$ be a normal unit vector to $\nil'$.
  Without loss of generality we can assume that $U\in
  \mathcal B$. We set $\Delta_0=-\sum_{V\in \mathcal
    B\setminus\{U\}} V^2$; we have $\Delta_0\in \mathfrak U
  (\nil')$ and $\Delta=-U^2 +\Delta_0$. Clearly
  $[\Delta_0]\le 2(k -1)$ and $[\Delta_0^r]\le 2(k -1)r$. It
  follows by~\pref{eq:bound} that there exists a constant
  $C:= C(k,r,\alpha)>0$ such that
  \begin{multline}
    \label{eq:bound2}
    \qquad \quad \Vert \pi_*(\Delta_0^r) G_{X}f \Vert\le
    C \,  \max  \{1, \delta_{\mathcal O}^{-1-2(k -1)r}\} \cdot \\
    \cdot\sum_{m=0}^{2(k -1)r} \Vert
    \left(I-\pi(\Delta_0)\right)^{\frac{m+\alpha}{2}} \pi_*(
    \delta^{m} \Delta_0^r) f\Vert. \qquad \quad
  \end{multline}
  
  Since $U$ is a unit vector orthogonal to $\nil'$, there
  exists $W\in \nil'$ such that
  \begin{equation}
    U= \<X,U\>^{-1} (X-W)\,.
  \end{equation}
  Since $\nil'\subset\nil$ is an ideal and $\pi(X)G_Xf=f$,
  it follows by estimate \pref{eq:Xproj} in Lemma
  \ref{lemma:cohomeqone} that there exists a constant
  $C':=C'(k,r,\alpha,X)>0$ such that
  \begin{equation}
    \label{eq:bound3}
    \begin{aligned}
      \Vert \pi_*(U)^{2r} G_{X}f \Vert & \leq \left(
        \frac{w_k(\Cal O)}{\delta_\Cal O} \right)^{2r}
      \Vert \pi_*(X-W)^{2r}  G_{X}f  \Vert \\
      &\leq C'\, \left( \frac{w_k(\Cal O)} {\delta_\Cal O}
      \right)^{2r} \left(\, \Vert f \Vert_{2r-1} + \Vert
        \pi_*(\Delta_0)^r G_{X}f \Vert\, \right)\,.
    \end{aligned}
  \end{equation}
  By bounds \pref{eq:bound2} and \pref{eq:bound3}, there
  exists $C'':=C''(k,r,\alpha,X)>0$ such that
  \begin{equation}
    \label{eq:bound4}
    \Vert  G_{X}f \Vert_{2r}  \leq C''\, \max\{1,w_k(\Cal
    O)^{2r}\}   
    \max\{ 1, \delta_{\Cal O} ^{-1-2rk} \}  
    \Vert f \Vert_{2kr +\alpha}\,.
  \end{equation}
  Finally by interpolation for all ~$ \alpha>\beta>0 $ such
  that $ \alpha> 1+k\beta $ there exists a constant
  $C_{\alpha,\beta}:=C_{\alpha,\beta}(k,X)>0$ such that
  \begin{equation}
    \label{eq:bound5}
    \Vert G_{X}f \Vert_{\beta} \leq C_{\alpha,\beta}\,
    \max\{1,w_k(\Cal O)^{\beta}\} \,
    \max \{1, \delta^{-1-k\beta}_{\Cal O}\} \, \Vert f\Vert_{\alpha}\,\,.
  \end{equation}
\end{proof}

\section{The cohomological equation for Nilflows}
\label{sec:flows}
\def\manif{\mathcal \Nil}
\subsection{The Howe-Richardson multiplicity formula}
\label{ssec:howe-richardson}

Let $\Nil$ be a $k$-step nil\-potent Lie group with a
minimal set of generators $ \{E_{1}, \dots,
E_{n}\}\subset \nil_{1} $.  Let $\Gamma$ be a lattice in
$\Nil$. The nilmanifold $\manif=\Gamma\backslash N$ carries
a unique $\Nil$-invariant probability measure $\mu$, locally
given by the Haar measure of $\Nil$; for each $X\in \nil$ we
denote by $(\phi^t_X)_{t\in \R}$ the the flow on~$\manif$
given by the right action of the one-parametre subgroup
$\{\exp t X \,\vert \,t\in \R\}$.

The Hilbert space $L^2(\mathcal N,\mu)$ decomposes under the
right action of $\Nil$ into a countable Hilbert sum
$\bigoplus_{i \in \N} H_i$ of irreducible closed subspaces
$H_i$.

\begin{sloppypar}
The Howe-Richardson multiplicity formula
\cite{MR0281842,MR0284546} tells us which irreducible
representations from the unitary dual $\hat{\Nil}$ of $\Nil$
appear in the decomposition $\bigoplus_{i \in \N} H_i$. To
state the formula, recall that a couple $(\chi, M)$ is
called a {\em maximal character} if there is $\lambda \in
\nil^*$ and a polarizing subalgebra $\mathfrak m \subset
\nil$ for $\lambda$ such that:
\end{sloppypar}
\begin{enumerate}
\item $M=\exp \mathfrak m$;
\item $\chi$ is the one-dimensional representation of $M$
  defined by $$\exp W\in M \mapsto e^{2 \pi \imath \lambda
    (W)};$$
\end{enumerate}
a maximal character $(\chi, M)$ is called {\em maximal 
  for $\Gamma$} if, in addition, the following conditions
are satisfied:
\begin{enumerate}
\item[(3)] $M$ intersects $\Gamma$ into a lattice, i.e.
  $M\cap \Gamma\backslash M$ is compact;
\item[(4)] $\chi$ is trivial on $M\cap \Gamma$: $\chi(M\cap
  \Gamma)=\{1\}$.
\end{enumerate}
As we recalled in \S\ref{ssec:kirillov}, when $(\chi, M)$ ranges
over the set of maximal characters, the family of
representations $\text{Ind}_M(\chi)$ exaust the unitary dual
$\hat \Nil$. The multiplicity formula states that
$\text{Ind}_M(\chi)$ appears in $L^2(\manif, \mu)$ iff
$(\chi, M)$ is a maximal integral character for $\Gamma$,
with multiplicity given by the cardinal of closed orbits
$xM$ of $M$ on $\manif$ such that $\chi( \text{Stab}_M(x) )
=\{1\}$. This discussion motivates the following:

\begin{definition} 
  A linear form $\lambda \in \nil^*$ is called \/ {\em
    integral (with respect to $\Gamma$)} if there exists a
  polarizing subalgebra $\mathfrak m \subset \nil$ such that
  the pair $\left(\exp (2\pi\imath \lambda), \exp \mathfrak
    m \right)$ is a maximal integral character for $\Gamma$.
  A coadjoint orbit $\Cal O \subset \nil$ will be called \/
  {\em integral (with respect to $\Gamma$)} if some
  $\lambda\in \Cal O$ is integral (hence all $\lambda\in
  \Cal O$ are).
\end{definition}

\begin{definition} 
  Let $\Gamma$ be a lattice in $\Nil$. A linear form
  $\lambda \in \nil^*$ is called \/ {\em weakly integral
    (with respect to $\Gamma$)} if $\lambda (E) \in \Z $ for
  all $E\in \log Z(\Gamma)$.  A coadjoint orbit $\Cal O
  \subset \nil$ will be called \/ {\em weakly integral (with
    respect to $\Gamma$)} if some $\lambda\in \Cal O$ is
  weakly integral (hence all $\lambda\in \Cal O$ are).
\end{definition}

Since $Z(\Gamma)= \Gamma \cap Z(\Nil) \subset \exp \mathfrak
m$ for any maximal polarizing subalgebra $\mathfrak m
\subset \nil$, any integral linear form (coadjoint orbit) is
also weakly integral. In fact, since the pair $\left(\exp
  (2\pi\imath \lambda), \exp \mathfrak m \right)$ is a
maximal integral character, $\exp\left(2\pi\imath
  \lambda(W)\right)=1$, hence $\lambda(W)\in \Z$, for all $W
\in \log Z(\Gamma)$.

It follows from the above discussion that all irreducible
unitary sub-representations occurring in $L^2(\manif, \mu)$
correspond are induced by Kirillov's by integral coadjoint
orbits.

\subsection{Diophantine elements}

Let $A$ the abelian group of rank $n$ given by
$A=\Nil/(\Nil,\Nil)$ and let $\bar \Gamma=
\Gamma/(\Gamma,\Gamma)$; then $\bar \Gamma$ is a lattice in
$A$.  Let $\{ \exp \bar E_1, \dots, \exp \bar E_n\}$ denote a
set of generators of $\bar \Gamma$, with $\bar E_i \in
\text{Lie}(A)\approx\nil/[\nil,\nil]$. It is plain that the
$\bar E_i$'s form a basis of $\nil/[\nil,\nil]$.

We shall say that an element $X\in\nil $ is {\em Diophantine (with
  respect to $ \Gamma $) of exponent~$ \tau \geq 0 $ } if
the projection $\bar X$ of $X$ in $\nil/[\nil,\nil]$ is
Diophantine of exponent $ \tau\geq 0 $ for the lattice $\bar
\Gamma$ in the standard sense: that is, if we let $\bar X=
\omega_1(X) \bar E_1 + \dots \omega_n(X) \bar E_n$, there
exists a constant $ K>0 $ such that, for all $
M=(m_{1},\dots, m_{n}) \in \Z^{n}\setminus \{0\} $,
\begin{equation}
  \label{eq:DC} \vert \<M,\Omega_{X}\> \vert := \vert \sum_{i=1}^
  {n} m_{i} \omega_{i}(X) \vert \,\, \geq \,\, \frac{K}{ \vert
    M\vert ^{n-1+\tau}}\,\,.
\end{equation}
The set of Diophantine elements $ X\in \nil $ of exponent $
\tau\geq 0 $ will be denoted by $ \text{DC}_{\tau}(\nil) $.
The subset of all Diophantine elements $ \text{DC}(\nil):=
\cup_{\tau} \text{DC}_{\tau}(\nil) \subset \nil $ has full
measure.
 
\begin{sloppypar}
  Let $E_1^1, E_2^1, \dots E_{n_1}^1, E_1^2, \dots,
  E_{n_2}^2, \dots, E_1^k, \dots, E_{n_k}^k$, (with
  $n_1=n$), a Malcev basis for $\nil$ through the descending
  central series $\nil_j$ and strongly based at $\Gamma$,
  that is a basis of $\nil$ satisfying the following
  properties:
  \begin{enumerate}
  \item if we drop the first $\ell$ elements of the basis we
    obtain a basis if a subalgebra of codimension $\ell$ of
    $\nil$.
  \item if we set $\mathcal E^j:=\{ E_1^j, \dots,
    E_{n_j}^j\}$ the elements of the set $\mathcal
    E^j\cup\mathcal E^{j+1} \cup \dots \cup\mathcal E^k $
    form a basis of $\nil_j$
  \item every element of $\Gamma$ can be written as a
    product
    $$
    \exp m_1^1 E_1^1 \dots \exp m_{n_1}^1 E_{n_1}^1\dots
    \exp m_1^k E_1^k \dots \exp m_{n_k}^k E_{n_k}^k
    $$
    with integral coefficients $m_i^j$.
  \end{enumerate}
  The existence of such a basis is obtained combining the
  proof of Theorems~1.1.13 and~5.1.6 of~\cite{CG:nilrep}.
\end{sloppypar}
It is also clear that we can assume that the basis $\{\bar
E_1, \dots, \bar E_n\}$ of $\nil/[\nil,\nil]$, appearing
above, is given by the projections in $\nil/[\nil,\nil]$ of
$E_1^1, E_2^1, \dots E_{n}^1$.

We define an Euclidean product on $\nil$ by making the basis
$E_j^k$ orthonormal. The norm of $\nil$ will be the one
induced by this product.

\begin{lemma}
  \label{lemma:cohomeqtwo} Let $ X\in
  \text{DC}_{\tau}(\nil)$.  Let $ \Cal O $ be a coadjoint
  orbit of maximal rank weakly integral (with respect to
  $\Gamma$) and let $ \pi \in \Pi_{\Cal O} $.  There exists
  a constant $ C_\Gamma>0 $ such that for all $ f\in
  W^{n-1+\tau} (H_{\pi}) $ we have:
  \begin{equation}
    \delta_{\Cal O}(X)^{-1} \Vert f \Vert \leq C_\Gamma \Vert f
    \Vert _{n-1+\tau} \,\,.
  \end{equation}
\end{lemma}

\begin{proof}
  Let $ \lambda \in \Cal O $.  Since $ \Cal O $ is weakly
  integral (with respect to $\Gamma$) and has maximal rank, we claim 
  that there exist $ Y\in {\Cal E}^{k-1} $ such that
  \begin{equation}
    M_{Y}:= \left( B_{\lambda}(E_{1}^1,Y), \dots, B_{\lambda}(E_
      {n}^1,Y) \right) \in \Z^{n}\setminus \{0\}\,\,.
  \end{equation}
  In fact, if $Y\in {\Cal E}^{k-1} $ then $ [E_{j}^1,Y]\in \nil_k $, since 
  ${\Cal E}^{k-1}\subset \nil_{k-1}$. By the
  Baker-Campbell-Hausdorff formula, we have $(\exp E_{j}^1,
  \exp Y) = \exp [E_{j}^1,Y] $. It follows that $[E_{j}^1,Y]\in \log
  Z(\Gamma)$ and, since $\lambda $ is weakly integral, we have
  $B_{\lambda}(E_ {j}^1,Y) =\lambda( [E_{j}^1,Y])\in \Z $
  for all $ j\in \{1,\dots,n\} $. If $
  B_{\lambda}(E_{j}^1,Y)=0 $, for all~$ Y\in {\Cal E}^{k-1}
  $ and all~$ j \in \{ 1,\dots ,n \} $, then $
  \nil_{k-1}^{\perp}(\Cal O) = \nil $ (see
  Lemma~\ref{lem:trivial}), contradicting the hypothesis
  that $ \lambda$ has maximal rank. The claim is proved.

  Since $ Y\in {\Cal E}^{k-1} $, by definition we have $
  \Vert Y\Vert =1 $.  Hence
  \begin{equation}
    \delta_{\Cal O}(X) \geq \vert B_{\lambda}(X,Y) \vert
    \geq \vert \<M_{Y},\Omega_{X}\> \vert \geq \frac{K}{%
      \vert M_{Y} \vert ^{n-1+\tau}}\,\,.
  \end{equation}

  Let $ \pi\in \Pi_{\Cal O} $ be an irreducible
  representation of the Lie algebra $ \nil $ on the Hilbert
  space $H:=H_{\pi} $.  Let~$ P = -
  (2\pi)^{-2}\sum_{i=1}^{n} [E_{i},Y] ^{2}\in {\mathfrak
    U}(\nil) $.  Since $ [E_{i},Y]\in Z(\nil) $, for all $
  i\in \{1,\dots,n\} $, the operator $ \pi_*(P) = \vert
  M_{Y}\vert ^{2} \I_{H} $.  Hence for any $ f \in H $ we
  have
  \begin{equation}
    \delta_{\Cal O}(X)^{-1} \Vert f \Vert \leq K^{-1} \Vert
    \vert M_{Y} \vert ^{n-1+\tau} f \Vert = K^{-1} \Vert \pi
    (P) ^{\frac{n-1+\tau}{2}} f \Vert \,\,.
  \end{equation}
  Since $ P $ has order $ 2 $ as an element of the
  enveloping algebra $\mathcal U(\nil)$ (in fact $ P\leq
  (2\pi)^{-2} \triangle $), if $ f\in W^{n-1+\tau}(H) $, the
  inequality $ \Vert \pi_*(P) ^{\frac{n-1+\tau}{2}} f \Vert
  \leq \Vert f \Vert _{n-1+\tau} $ holds.
\end{proof}

\subsection {Uniform estimates}

We prove that, for any Diophantine element, appropriate
Sobolev norms of the Green operators for the cohomological
equation, which have been constructed in each irreducible
unitary representation, are {\it bounded uniformly }over all
irreducible unitary representations. This step is crucial in
order to solve the cohomological equation on a nilmanifold
$\manif= \Gamma\backslash \Nil$ by ''gluing'' together the
solutions constructed in every irreducible
sub-representation of the representation of $\Nil$ on
$L^2(\manif,\mu)$.

\begin{lemma} \label{lemma:cohomeqthree} Let $ \Cal O $ be any
  coadjoint orbit and let $\pi \in \Pi_{\Cal O}$ be a
  unitary representation of $\Nil$.  For all $f\in
  W^1(H_{\pi})$, the following inequality holds:
  \begin{equation}
    w_Z(\Cal O) \Vert f \Vert \leq  \frac{1}{2\pi}\,\Vert f \Vert_1\,.
  \end{equation}
\end{lemma}
\begin{proof} By the definition \pref{eq:ws} of the
  non-negative number $w_Z(\Cal O)$, for each coadjoint
  orbit $\Cal O$ of $\Nil$, there exists an element $Z\in
  Z(\nil)$ such that $ \Vert Z \Vert =1$, with respect to
  the fixed euclidean norm on $\nil$, and $\vert \lambda(Z)
  \vert = w_Z(\Cal O)$, for any $\lambda \in \Cal O$. In
  addition, since $Z\in Z(\nil)$, for any $\pi \in \Pi_{\Cal
    O}$ we have that $\pi(Z) =2\pi \imath \lambda(Z)$. It
  follows that, for any $f\in W^1(H_{\pi})$,
  \begin{equation}
    w_Z(\Cal O) \Vert f \Vert =  
    \frac{1}{2\pi}\,  \Vert \pi_*(Z) f  \Vert \leq   \frac{1}{2\pi}\,
    \Vert f \Vert_1\,,
  \end{equation}
  as claimed.
\end{proof}

\begin{corollary}
  \label{cor:cohomeqrep}
  Let $ X\in \text{DC}_{\tau}(\nil) $ .  Let $ \Cal O $ be a
  non-trivial coadjoint orbit , integral with respect to a
  lattice $\Gamma$, and let $\pi \in \Pi_{\Cal O} $ an
  irreducible representation of the Lie algebra $ \nil $ on
  a Hilbert space $ H_{\pi} $.

  \begin{enumerate}

  \item If $ \alpha>n + \tau -1/2$ and $ \beta<-1/2 $, there
    exists $ C'_\Gamma:=C'_\Gamma(X,\alpha)>0 $ such that, for all $ f\in
    W^{\alpha}(H_{\pi}) $,
    \begin{equation*}
      \Vert G_{X}f \Vert _{\beta} \le C'_\Gamma \, \Vert f \Vert _
      {\alpha} \,\,;
    \end{equation*}
  \item if $ f\in W^{\alpha}(H_{\pi}) $, $ \alpha>n+\tau $,
    and $D(f)=0 $, for all $D\in {\Cal I}^
    {\alpha}_{X}(H_{\pi}) $, then $ G_{X}f \in W^{\beta}(H_
    {\pi}) $, for all $ \beta<[\alpha-(n+\tau)]
    [(n+\tau)k+1]^ {-1} $ and there exists $
    C''_\Gamma:=C''_\Gamma(X,\alpha,\beta)>0 $ such that
    \begin{equation*}
      \Vert G_{X}f \Vert _{\beta} \le C''_\Gamma\, \Vert f \Vert _
      {\alpha} \,\,.
    \end{equation*}
  \end{enumerate}
\end{corollary}
\begin{proof}
  The proof is by induction on the degree of nilpotency
  $k\in\N\setminus\{0\}$ of the nilpotent group $\Nil$. If
  $k=1$ the group $\Nil$ is abelian, all irreducible unitary
  representations are one-dimensional and (weakly) integral.
  The Diophantine condition immediately implies that the
  statement holds in this case with the exception of the
  trivial representation. Let us assume that the statement
  holds for all lattices in any $(k-1)$-step nilpotent group
  and let $\Nil$ be a $k$-step nilpotent group. If the
  coadjoint orbit $\Cal O$ is of maximal rank and integral
  with respect to a lattice $\Gamma\subset \Nil$, then the
  statement follows from Theorem \ref{th:cohomeq},
  Lemma~\ref{lemma:cohomeqtwo} and
  Lemma~\ref{lemma:cohomeqthree}.

  If $\Cal O$ is integral (with respect to $\Gamma$), but it
  is not of maximal rank, then for any $\lambda\in \Cal O$,
  by Lemma~\ref{lem:nonmaxrk} the restriction $\lambda\vert
  \nil_k$ is identically zero and any irreducible
  representation $\pi\in \Pi_{\Cal O}$ factors through an
  irreducible representation $\pi'$ of the $(k-1)$-step
  nilpotent group $\Nil':=\Nil/\exp\nil_k$, induced by the
  linear functional $\lambda' \in ( \nil/\nil_k )^*$.  Since
  $\lambda\in \nil^*$ is integral (with respect to
  $\Gamma$), there exists a polarizing subalgebra $\mathfrak
  m$ such that $\left(\exp(2\pi\imath \lambda), \exp
    \mathfrak m\right)$ is a maximal integral character for
  $\Gamma$. Since $\lambda\vert \nil_k \equiv 0$, the
  subalgebra $\mathfrak m':= \mathfrak m/\nil_k$ is
  polarizing for $\lambda'$ and $\left(\exp(2\pi\imath
    \lambda'), \exp \mathfrak m' \right)$ is a maximal
  integral character for the lattice $\Gamma':= \Gamma/
  \Gamma\cap \exp\nil_k$. By the induction hypothesis, the
  statement holds also in this case.
\end{proof}

\subsection {Proof of Theorems \ref{th:invardistrib} and \ref{th:main}} 
Let $L^2(\manif, \mu)= \oplus_{i\in \N} H_i$ be the
decomposition of the space of square integrable functions on
$\manif$ into irreducible components under the right action
of $\Nil$.  Such a decomposition induces a decomposition of
the subspace of smooth vectors,
$$
C^{\infty}(\manif) = \bigoplus_{i\in \N} C^{\infty}(H_i)\,,
$$
and of the Sobolev spaces: for each $\alpha\in \R$, there is
an {\it orthogonal }splitting
\begin{equation}
\label{eq:Sobsplit}
W^\alpha(\manif) = \bigoplus _{i\in \N} W^\alpha(H_i)\,.
\end{equation}
Consequently, for any $X\in \nil$, the spaces $\Cal I_X
(\manif) \subset \Cal D' (\manif)$ of all $X$-invariant
distributions and the subspace $\Cal I^{\alpha}_X (\manif):=
\Cal I_X (\manif)\cap W^{ -\alpha}(\manif)$ split as
follows:
\begin{equation}
\label{eq:invdistdec}
\begin{aligned}
\Cal I_X (\manif)&= \bigoplus_{i\in \N} \Cal I_X (H_i) \,; \\
\Cal I^\alpha_X (\manif)&= \bigoplus_{i\in \N} \Cal I^\alpha_X (H_i) \,. 
\end{aligned}
\end{equation}
Let $ X $ be irrational. By Lemma~\ref{lem:invdist1} and by
the decompositions~\pref{eq:invdistdec}, the space $\Cal I_X
(\manif)$ admits a countable basis $\Cal B_X (\manif)
\subset \Cal I^\alpha_X$, for any $\alpha >1/2$. This proves
Theorem~\ref{th:invardistrib}

\smallskip Let $ X\in \text{DC}_{\tau}(\nil) $. Then the
Green operator for the cohomological equation on $\manif$
can be constructed as follows.  Let $\alpha>n+\tau -1/2$ and
$\beta < -1/2$. For any $i\in \N$, let $G_X^{(i)}$ be the
Green operator for the cohomological equation for $X$ in the
irreducible representation $H_i$. By
Corollary~\ref{cor:cohomeqrep}, the operators $G_X^{(i)}:
W^\alpha(H_i) \to W^\beta(H_i)$ are bounded and have
uniformly bounded norms, in the sense that there exists
$C'_\Gamma>0$ such that $\Vert G_X^{(i)}\Vert
_{\alpha,\beta} \leq C'_\Gamma$.

For each $i\in \N$ and $\alpha\in \R$, let $p_\alpha^{(i)}:
W^\alpha(\manif) \to W^\alpha(H_i)$ be the orthogonal
projection and $j_\alpha^{(i)}: W^\alpha(H_i) \to
W^\alpha(\manif)$ the embedding, determined by the
orthogonal splitting~\pref{eq:Sobsplit}. The operator $G_X:
W^\alpha(\manif) \to W^\beta(\manif)$ defined as
$$
G_X : = \bigoplus_{i\in \N}    
j_\beta^{(i)} \circ G_X^{(i)} \circ p_\alpha^{(i)}
$$
is a well-defined, bounded Green operator for the
cohomological equation for the Diophantine vector field
$X\in \text{DC}_{\tau}(\nil) $ on $\manif$. In fact, by the
orthogonality of the splitting~\pref{eq:Sobsplit},
$$
\Vert G_X \Vert^2_{\alpha,\beta} \leq \sum_{i\in \N} \Vert
j_\beta^{(i)} \circ G_X^{(i)} \circ p_\alpha^{(i)} \Vert^2
_{\alpha,\beta} \leq \max_{i\in\N} \Vert G_X^{(i)}
\Vert_{\alpha,\beta}^2 \leq (C'_\Gamma)^2\,.
$$
Finally, if $\alpha>n+\tau$ and $f\in W^{\alpha}(\manif)$
belongs to the kernel of all $X$-invariant distributions $D
\in \Cal I_X^\alpha(\manif)$, then $u:=G_X(f) \in
W^\beta(\manif)$ for any $ \beta<[\alpha-(n+\tau)]
[(n+\tau)k+1]^ {-1} $. By construction,
$$
u = \sum_{i\in \N} j_\beta^{(i)} \circ G_X^{(i)} \left(
  p_\alpha ^{(i)} (f) \right)\,.
$$
By the splitting~\pref{eq:invdist}, it follows that $D_i
(p_\alpha^{(i)}(f))=0$, for all $D_i\in \Cal
I_X^\alpha(H_i)$ and for each $i\in \N$, since $D(f)=0$ for
all $D\in \Cal I_X^\alpha(\manif)$. Hence, by
Corollary~\ref{cor:cohomeqrep}, there exists a constant
$C''_\Gamma>0$ such that, for all $i\in \N$,
$$
\Vert G_X^{(i)} \left( p_\alpha^{(i)}(f) \right)
\Vert_{\beta} \leq C''_\Gamma \,\Vert
p_\alpha^{(i)}(f)\Vert_{\alpha}\,.
$$
Again by the orthogonality of the
splitting~\pref{eq:Sobsplit}, we have
$$
\Vert u\Vert^2_\beta = \sum_{i\in \N} \Vert G_X^{(i)} \left(
  p_\alpha^{(i)}(f)\right)\Vert^2_\beta \leq (C''_\Gamma)^2
\sum_{i\in \N} \Vert p_\alpha^{(i)}(f) \Vert^2_\alpha \leq
(C''_\Gamma)^2 \Vert f\Vert^2_\alpha\,.
$$
The proof Theorem \ref{th:main} is concluded.

\vskip 1cm
\bibliography{biblio}
\bibliographystyle{amsalpha}

\end{document}

%% file: env_def.tex
\newtheorem{theorem}{Theorem}[section]
\newtheorem{lemma}[theorem]{Lemma}
\newtheorem{corollary}[theorem]{Corollary}

\newtheorem{question}{Question}
\newtheorem{conjecture}{Conjecture}
\newtheorem{definition}[theorem]{Definition}

\theoremstyle{definition}

%% file: symb_def.tex
\newcommand{\field}[1]{\mathbb{#1}}
\newcommand{\R}{\field{R}}
\newcommand{\N}{\field{N}}
\newcommand{\C}{\field{C}}
\newcommand{\Z}{\field{Z}}
\newcommand{\Q}{\field{Q}}
\newcommand{\T}{\field{T}}
\newcommand{\Cal}{\mathcal}
\newcommand{\manif}{\Cal N}

\newcommand{\ad}{\hbox{ad}}
\newcommand{\Ad}{\hbox{Ad}}

\newcommand{\I}{\text{\rm Id}}
\newcommand{\<}{{\langle}}
\newcommand{\m}{\mathfrak{m}}

\newcommand{\nil}{\mathfrak{n}}
\newcommand{\Nil}{{N}}
\newcommand{\pref}[1]{(\ref{#1})}

\renewcommand{\>}{{\rangle}}
\renewcommand{\|}{\,\Vert\,}
